\documentclass{ACmod}

\DeclareFontFamily{U}{rsf}{}
\DeclareFontShape{U}{rsf}{m}{n}{
  <5> <6> rsfs5 <7> <8> <9> rsfs7 <10-> rsfs10}{}
\DeclareMathAlphabet{\mathscr}{U}{rsf}{m}{n}

\makeatletter
\def\operator@font{\sf}
\makeatother

\newcommand{\cO}{{\mathscr O}}

\newcommand{\deform}{\mathsf{def}}

\DeclareMathOperator{\HH}{HH}
\DeclareMathOperator{\HT}{HT}

\DeclareMathOperator{\HOmega}{H}

\DeclareMathOperator{\Ext}{Ext}

\DeclareMathOperator{\codim}{codim}

\renewcommand{\phi}{\varphi}

\usepackage{tikz-cd}
\usepackage{amsmath}

\usepackage{xypic}
\usepackage{tikz}
\usepackage[tocflat]{tocstyle}
\usetikzlibrary{positioning}
\usepackage{amscd}
\usetikzlibrary{arrows,shapes,positioning}
\usetikzlibrary{decorations.markings}
\tikzstyle arrowstyle=[scale=1]
\tikzstyle directed=[postaction={decorate,decoration={markings,
    mark=at position .65 with {\arrow[arrowstyle]{stealth}}}}]
\tikzstyle reverse directed=[postaction={decorate,decoration={markings,
    mark=at position .65 with {\arrowreversed[arrowstyle]{stealth};}}}]
\usepackage[tocflat]{tocstyle}

\tikzset{cross/.style={cross out, draw=black, minimum size=2*(#1-\pgflinewidth), inner sep=0pt, outer sep=0pt},
cross/.default={1pt}}

\begin{document}

\title{The orbifold Hochschild product for Fermat hypersurface}

\author[Shengyuan Huang and Kai Xu]{%
Shengyuan Huang and Kai Xu}

\address{
\textsc{Shengyuan Huang: School of Mathematics,
University of Birmingham
\newline
Edgbaston, Birmingham, B15 2TT, UK}
\newline
Email address: \texttt{s.huang.5@bham.ac.uk}
\newline
\textsc{Kai Xu: Mathematics Department,
Harvard University
\newline
1 Oxford St, Cambridge, MA 02138, US}
\newline
Email address: \texttt{kaixu@math.harvard.edu}
}

\begin{abstract}
{\sc Abstract:}
Let $G$ be an abelian group acting on a smooth algebraic variety $X$. We investigate the product structure and the bigrading on the cohomology of polyvector fields on the orbifold $[X/G]$, as introduced by C\u{a}ld\u{a}raru and Huang. In this paper, we provide many new examples given by quotients of Fermat hypersurfaces, where the product is shown to be associative. This is expected due to the conjectural isomorphism at the level of algebras between the cohomology of polyvector fields and the Hochschild cohomology of orbifolds. We prove this conjecture for the Calabi-Yau Fermat hypersurface orbifold. We also show that for Calabi-Yau orbifolds, the multiplicative bigrading on the cohomology of polyvector fields agrees with what is expected in homological mirror symmetry. 
\end{abstract}

\maketitle
\setcounter{tocdepth}{1}
\tableofcontents
\vfill

\section{Introduction}

\paragraph Let $X$ be a smooth algebraic variety over a field of characteristic zero. The Hochschild cohomology of $X$ is well-understood due to the work of Swan, Kontsevich, Calaque-Van den Bergh, and many other mathematicians. The HKR map~\cite{S} is an isomorphism of vector spaces
$$\HT^*(X)\stackrel{\cong}{\rightarrow}\HH^*(X)$$
between the {\em polyvector field cohomology}
$$\HT^*(X)=\bigoplus_{p+q=*}H^p(X,\wedge^qT_X)$$
and the Hochschild cohomology $\HH^*(X)$.

The Hochschild cohomology $\HH^*(X)$ is by definition $\Ext^*(\Delta_*\cO_X,\Delta_*\cO_X)$, where $\Delta: X\hookrightarrow X\times X$ is the diagonal embedding. Therefore its classes can be composed using the Yoneda product. There is a wedge product on polyvector fields. However, the HKR map is not an isomorphism of algebras in general. Kontsevich~\cite{K} constructed a highly nontrivial map $\HT^*(X)\rightarrow\HH^*(X)$ which was proven to be an isomorphism of algebras by Calaque and Van den Bergh~\cite{CV}.

\paragraph
Mathematicians start to study the multiplicative structure of the Hochschild cohomology for orbifolds recently. Some progress has been made by Arinkin-C\u{a}ld\u{a}raru-Hablicsek~\cite{ACH}, Negron-Schedler~\cite{NS}, and C\u{a}ld\u{a}raru-Huang~\cite{CH}. This product structure is the mirror to the quantum cup product on the mirror hence it reveals deep enumerative geometry under mirror symmetry.

Let $G$ be a finite group acting on a smooth algebraic variety $X$ over a field of characteristic zero. The Hochschild cohomology $\HH^*([X/G])$ of the orbifold $[X/G]$ has a natural algebra structure. Having an explicit formula for the product structure of the algebra $\HH^*([X/G])$ would yield many possible applications in homological mirror symmetry and the crepant resolution conjecture.

\paragraph
Arinkin, C\u{a}ld\u{a}raru, and Hablicsek~\cite{ACH} gave an explicit
decomposition of the Hochschild cohomology of $[X/G]$ in terms of polyvector field cohomology for orbifolds. They showed that there exists a graded vector space isomorphism
\[
\HH^*([X/G]) \cong \HT^*([X/G])\stackrel{\deform}{=}\left (\bigoplus_{g\in
    G}\bigoplus_{p+q=*} H^{p-c_g}(X^g,\wedge^q
  T_{X^g}\otimes\omega_g)\right )^G\!\!\!\!\!,
\]
where $X^g$ is the fixed locus of $g\in G$, $c_g$ is the codimension
of $X^g$ in $X$, and $\omega_g$ is the dualizing sheaf of the
inclusion $X^g\hookrightarrow X$. This is now called the orbifold HKR isomorphism.

The Hochschild cohomology $\HH^*([X/G])$ has a natural product. When $G$ is abelian, C\u{a}ld\u{a}raru and Huang~\cite{CH} defined a product on
\[
 \HT^*(X;G)\stackrel{\deform}{=}\bigoplus_{g\in
    G}\bigoplus_{p+q=*} H^{p-c_g}(X^g,\wedge^q
  T_{X^g}\otimes\omega_g).
\]
Note that $\HT^*(X;G)$ carries a natural $G$-action and the $G$-invariant part is $\HT^*([X/G])$. This product on $\HT^*(X;G)$ is the wedge product on $\HT^*(X)$ in the case where $G$ is trivial. The authors in~\cite{CH} conjecture that the two algebras $\HH^*([X/G])$ and $\HT^*([X/G])$ are isomorphic in a highly nontrivial way which generalizes the isomorphism of algebras~\cite{K} in the case where $G$ is trivial.

\paragraph{\bf The associativity.} The first evidence that would be needed for such an isomorphism is that the product on $\HT^*(X;G)$ is associative. The authors of ~\cite{CH} showed that the product they defined is associative when the Bass-Quillen class vanishes. A few examples have been computed in~\cite{CH} and the Bass-Quillen class vanishes there, but the size of the group $G$ is small in those examples.

In this paper we consider the $(\mathbb{Z}/d\mathbb{Z})^{n-1}$ action on the degree $d$ Fermat hypersurface in $\mathbb{P}^n$. We prove that the product is associative in this case. This provides examples such that the product on $\HT^*(X; G)$ is associative with an arbitrarily large group $G$.

\paragraph{\bf Theorem A.}
{\em Let $[x_0:\cdots:x_n]$ be the homogenous coordinates on $\mathbb{P}^n$. The degree $d$ Fermat hypersurface $X$ in $\mathbb{P}^n$ is defined by $\sum_{j=0}^{n}x_j^{d}=0$. The group $G=(\mathbb{Z}/d\mathbb{Z})^{n-1}$ acts on $X$. Let $\zeta=\exp^{2\pi i/d}$ be the root of unity. An element of $G$ is of the form $g=(\zeta^{a_0},\zeta^{a_1},\cdots,\zeta^{a_{n-1}},1)$, where $a_j\in\mathbb{Z}/d\mathbb{Z}$ and $\sum_{j=0}^{n-1}a_j=0$. The group action of $G$ on $X$ is defined by
$$g\cdot[x_0:\cdots:x_n]=[\zeta^{a_0}x_0:\cdots:\zeta^{a_{n-1}}x_{n-1}:x_n].$$
Then the Bass-Quillen classes~\cite{CH, H} associated to the sequences $X^{g,h}\hookrightarrow X^{gh}\hookrightarrow X$ and $X^{g,h}\hookrightarrow X^{g}\hookrightarrow X$ vanish for all $g,h\in G$. Therefore the product~\cite{CH} defined on
\[
 \HT^*(X;G)=\bigoplus_{g\in
    G}\bigoplus_{p+q=*} H^{p-c_g}(X^g,\wedge^q
  T_{X^g}\otimes\omega_g)\!\!\!\!\!,
\]
is associative.
}

\paragraph
The product in~\cite{CH} is defined by abstract tools in derived algebraic geometry. It is difficult to compute the product explicitly in general. However, in Section 7 of {\em loc. cit.} the authors defined a new product which we call the {\em simplified product} on $\HT^*(X;G)$
by explicit formulas. We call the original product the {\em unsimplified product}. Conjecturally the simplified product is equal to the unsimplified product for Calabi-Yau orbifolds. Here Calabi-Yau orbifolds are the direct generalization of Calabi-Yau varieties:

\begin{Definition}
    A smooth proper orbifold (or Deligne-Mumford stack) $Y$ is said to be Calabi-Yau if its canonical bundle is trivial, i.e. $\omega_Y\simeq \mathcal{O}_Y$ as an orbifold line bundle. 
\end{Definition}

\paragraph
One may directly check that a global quotient orbifold $[X/G]$ is Calabi-Yau if and only if $X$ is Calabi-Yau and the $G$ action on $X$ preserves the holomorphic volume form $\Omega_X\in\Gamma(X,\omega_X)$. A first example is the Fermat quintic with the $(\mathbb{Z}/5\mathbb{Z})^3$ action. In this paper, we show that in this example the simplified and unsimplified products agree.

\paragraph{\bf Theorem B1.}
{\em In the case of Fermat quintic, i.e., $d=5$ and $n=4$, the unsimplified and the simplified products on $\HT^*(X;G)$ are equal.
}

For a general Calabi-Yau Fermat hypersurface, we are not able to prove the simplified and unsimplified product on $\HT^*(X;G)$ agree. However, we can prove that they agree after taking $G$-invariants.

\paragraph{\bf Theorem B2.}
{\em In the case of Calabi-Yau Fermat hypersurface, i.e., $d=n+1$, the unsimplified and the simplified products on $\HT^*([X/G])=\HT^*([X/G])^G$ are equal.
}

\paragraph{\bf The multiplicative bigrading.}\label{age}
To prove Theorem B2, we need to study the multiplicative bigrading on $\HT^*(X;G)$ and apply Theorem C below. The authors in~\cite{CH} defined a new bigrading on $\HT^*(X;G)$
as follows
\[
 \HT^{q,p}(X;G)=\bigoplus_{g\in
    G}H^{p-\iota(g)}(X^g,\wedge^{q-c_g+\iota(g)}
  T_{X^g}\otimes\omega_g)\!\!\!\!\!,
\]
where $\iota(g)$ is the age of $g\in G$~\cite{CR,FG}. The age $\iota(g)$ of a group element $g$ is defined as follows. Let X be an algebraic variety of dimension $D$ with the action of a finite group $G$. For $g\in G$ and a point $x\in X^g$, let $\lambda_1,\cdots,\lambda_D$ be the eigenvalues of the action of $g$ on the tangent space $T_{X,x}$; note that they are roots of unity. Write $\lambda_j=e^{2\pi i r_j}$ where $r_j$ is a rational number in the interval $[0,1)$. The age $\iota(g,y)$ of $g$ in $y$ is the
rational number $\sum r_j$. The age $\iota(g,y)$ only depends on the connected component $Z$ of $Y^g$ in which $y$ lies~\cite{FG}. For simplicity, we write $\iota(g)$ instead of $\iota(g,Z)$ when the connected component $Z$ is clear from the context.

%let $\lambda_j=e^{2\pi i r_j}$ be the eigenvalue of the action of $g$ on the tangent space where $r_j\in\mathbb{Q}\cap[0,1)$, then $\iota(g)=\sum r_j$. The simplified product preserves the $(q,p)$ bidegree.

When $[X/G]$ is Calabi-Yau, we show that the bigrading defined above agrees with the bigrading on the orbifold singular cohomology of the mirror in the sense of Theorem C below. To explain Theorem C, we start with the variety case first.

\begin{Definition}
For a smooth algebraic variety $X$, define
$$\HT^{q,p}(X)\stackrel{\deform}{=}H^{p}(X,\wedge^qT_X),$$
and
$$\HOmega^{q,p}(X)\stackrel{\deform}{=}H^p(X,\wedge^q\Omega_X).$$
\end{Definition}

For a Calabi-Yau variety $X$ of dimension $n$, the dimension of the cohomology of polyvector fields is related to the dimension of singular cohomology of $X$ because of the identification

$$\wedge^{q}T_X\cong\wedge^{n-q}\Omega_X.$$
It implies $\HT^{q,p}(X)\cong\HOmega^{n-q,p}(X)$.

\label{mirror}
Mirror symmetry predicts that $\HOmega^{n-q,p}(X)\cong\HOmega^{q,p}(\check{X})$ if the Calabi-Yau variety $X$ has a mirror $\check{X}$. Hence mirror symmetry expects the isomorphism $\HT^{q,p}(X)\cong\HOmega^{n-q,p}(X)\cong\HOmega^{q,p}(\check{X})$. In fact homological mirror symmetry predicts that
$$\HT^*(X)=\bigoplus_{p+q=*}\HT^{q,p}(X)$$
should be isomorphic to
$$H^*(\check{X},\mathbb{C})=\bigoplus_{p+q=*}\HOmega^{q,p}(\check{X})$$
as bigraded algebras~\cite{Kont, Se}, where the product on $H^*(\check{X},\mathbb{C})$ should be the quantum product rather than the singular cohomology product.

\paragraph For an orbifold $[X/G]$, there is an orbifold version of Hodge decomposition arising from orbifold singular cohomology

\[
H^{*}([X/G],\mathbb{C})=\left( \bigoplus_{g\in G}\bigoplus_{p+q=*}H^{p-\iota(g)}(X^g,\wedge^{q-\iota(g)}\Omega_{X^g}) \right)^G
\]
of $[X/G]$ defined by Chen and Ruan~\cite{CR}. This has a bigrading given as follows.

\begin{Definition}
For a global quotient orbifold $[X/G]$, define
$$\HOmega^{q,p}(X;G)=\bigoplus_{g\in G}H^{p-\iota(g)}(X^g,\wedge^{q-\iota(g)}\Omega_{X^g}),$$
\[  \HOmega^{q,p}([X/G])=\HOmega^{q,p}(X;G)^G,
\]
and
$$\HT^{q,p}(X;G)=\bigoplus_{g\in G}H^{p-\iota(g)}(X^g,\wedge^{q+\iota(g)-c_g}T_{X^g}\otimes\omega_g),$$
$$\HT^{q,p}([X/G])=\HT^{q,p}(X;G)^G.$$
\end{Definition}

Both of the bigradings are multiplicative, i.e., they are preserved by the corresponding product structures respectively. The dimension of $\HOmega^{q,p}([X/G])$ is the $(q,p)$-th orbifold Hodge number of $[X/G]$. When $[X/G]$ is Calabi-Yau, i.e. its canonical bundle is trivial as an orbifold line bundle, we prove that the bigrading for $\HT^*([X/G])$ coincides with the bigrading on the singular cohomology of its mirror.

\paragraph{\bf Theorem C.}
{\em Let $[X/G]$ be a Calabi-Yau orbifold of dimension $n$. Then $\HT^{q,p}(X;G)\cong\HOmega^{n-q,p}(X;G)$ and $\HT^{q,p}([X/G])\cong\HOmega^{n-q,p}([X/G])$.
}

If $[X/G]$ is Calabi-Yau of dimension $n$ and has a mirror $[Y/H]$, then $\HOmega^{n-q,p}([X/G])$ should be identified with $\HOmega^{q,p}([Y/H])$. Theorem C shows that
\[
\HT^*([X/G])=\bigoplus_{p+q=*}\HT^{q,p}([X/G])
\]
is identified with
\[ H^*([Y/H],\mathbb{C})=\bigoplus_{p+q=*} \HOmega^{q,p}([Y/H])
\]
as bigraded vector spaces. This provides evidence that the multiplicative bigrading we put on $\HT^*(X;G)$ is the correct one. Note that Theorem C only requires $[X/G]$ to be Calabi-Yau, not necessarily to be the Fermat hypersurface orbifold.

Theorems B2 and C above and Conjecture A in~\cite{CH} suggest that the Hochschild
cohomology of Calabi-Yau orbifolds should carry a
multiplicative bigrading.

\paragraph We return to the Calabi-Yau Fermat hypersurface orbifold case, where $X$ is the Calabi-Yau Fermat hypersurface of degree $d$ and $G$ is the group $(\mathbb{Z}/d\mathbb{Z})^{d-2}$. We compute the product on $\HT^*([X/G])$ explicitly.

For $g,h\in G$ let $\epsilon(g,h)=\iota(g)+\iota(h)-\iota(gh)$. We define a modified algebra structure on $\HT^*([X/G])$ by 
$$\alpha_g\circ \beta_h=(-1)^{\epsilon(g,h)}\alpha_g\cdot\beta_h.$$ 
Denote this new algebra by $(\HT^*([X,G]),\circ)$. Note that the sign is not surprising because Fantechi and G\"{o}ttsche~\cite[Definition 3.9]{FG} have also introduced a similar sign in their study of orbifold singular cohomology. The product is graded commutative after the sign is introduced.

%The quotient variety $X/G$ is a singular variety whose crepent resolution is the mirror quintic $\check{X}$. Then we have the following identification of algebras

%$$\HH^*([X/G])\cong\HH^*(\check{X})\cong H_{FJRW}^*(W),$$
%where $H_{FJRW}^*(W)$ is the state space of the FJRW theory of the Fermat polynomial $W=\sum_{i=0}^4 x_i^5$. The first isomorphism is due to Bridgeland, King, and Reid~\cite{BKR}. The second isomorphism is due to mirror symmetry and more details are explained in the last section of this paper. We obtain the following corollary of Theorem D by direct computation.
We study this modified product structure in detail and obtain Theorem D below.
%\paragraph{\bf Corollary A.}
%{The orbifold polyvector field $\HT^*([X/G])$ of Fermat quintic orbifold is isomorphic to the state space of the FJRW theory $H^*_{FJRW}(W)$ of the Fermat polynomial $W=\sum_{i=0}^4 x_i^5$ as algebras.}

\paragraph{\bf Theorem D.}
{\em In the case of Calabi-Yau Fermat hypersurface orbifold, the orbifold polyvector field with the modified product $(\HT^*([X/G]),\circ)$ is isomorphic to the orbifold Hochschild cohomology $\HH^*([X/G])$ as algebras.
}
 
Theorem D provides a positive answer to Conjecture A in~\cite{CH} in the case of the Calabi-Yau Fermat hypersurface orbifold. We discuss connections in the field of mirror symmetry at the last of this paper.

\paragraph {\bf Plan of the paper.}
In Section 2 we compute the fixed locus of an element $g\in G$. Then we prove Theorem A.

In Section 3 we study the group $G=(\mathbb{Z}/5\mathbb{Z})^3$. We classify the elements of the group into four types. The result will be used in Sections 4-5.

In Section 4 we study the simplified product in~\cite{CH}. The definition of the simplified product depends on a cohomology class which is introduced in~\cite{FG}. We show that the simplified product can be simplified further when the class is trivial. We use the results to study the simplified product in the case of the Fermat quintic orbifold.

In Section 5 we prove Theorem B1.

In Section 6 we prove Theorem C.

In section 7 we prove Theorems B2 and D.

\paragraph {\bf Acknowledgments.} We would like to thank Andrei C\u{a}ld\u{a}raru
for introducing the subject of orbifold Hochschild cohomology and thank Tyler Kelly for his useful suggestions on writing. We are also grateful to the referee for careful reading and many suggestions and comments which greatly improved this paper.

The first author was partially supported by the UKRI Future Leaders Fellowship through grant number MR/T01783X/1.

\section{Proof of Theorem A}
We study the fixed locus of an element $g\in G$. Then we prove Theorem A.

\paragraph Let $[x_0:\cdots:x_n]$ be the homogenous coordinates on $\mathbb{P}^n$. The degree $d$ hypersurface $X$ in $\mathbb{P}^n$ is defined by $\sum_{j=0}^{n}x_j^{d}=0$. The group $G=(\mathbb{Z}/d\mathbb{Z})^{n-1}$ acts on $X$. Let $\zeta=\exp^{2\pi i/d}$ be the root of unity. An element of $G$ is of the form $g=(\zeta^{a_0},\zeta^{a_1},\cdots,\zeta^{a_{n-1}},1)$, where $a_j\in\mathbb{Z}/d\mathbb{Z}$ and $\sum_{j=0}^{n-1}a_j=0$. The group action is defined by
$$g\cdot[x_0:\cdots:x_n]=[\zeta^{a_0}x_0:\cdots:\zeta^{a_{n-1}}x_{n-1}:x_n].$$
We need to study the fixed locus of $g\in G$ before we study the product on $\HT^*(X;G)$.

\paragraph{\bf The fixed locus.}\label{fixed locus}
The fixed locus depends on the numbers of $a_0$, ..., $a_{n-1}$ that are equal to each other. We compute an example and one can generalize the proof to the general case. For example, choose $d=5$, $n=8$, and $g=(\zeta,\zeta,\zeta,\zeta^2,\zeta^2,\zeta^2,\zeta^2,\zeta^{-1},1)$. Let $[x_0:\cdots :x_8]$ be a point in the fixed locus of $g$. Then
$$g\cdot[x_0:\cdots:x_8]=[\zeta x_0:\zeta x_1:\zeta x_2:\zeta^2 x_3:\zeta^2 x_4:\zeta^2:x_5:\zeta^2 x_6: \zeta^{-1}x_7:x_8].$$

By the definition of the homogenous coordinates, there exists a nonzero number $\lambda$ such that
$$(\zeta x_0,\zeta x_1,\zeta x_2,\zeta^2 x_3,\zeta^2 x_4,\zeta^2x_5,\zeta^2 x_6, \zeta^{-1}x_7,x_8)=(\lambda x_0,\cdots,\lambda x_8).$$

If any of $x_0$, $x_1$, and $x_2$ is nonzero, then $\lambda$ has to be $\zeta$. We can conclude $x_3=\cdots=x_8=0$.

If any of $x_3$, $x_4$, $x_5$, and $x_6$ is nonzero, then $\lambda$ has to be $\zeta^2$. We can conclude $x_0=x_1=x_2=0$ and $x_7=x_8=0$.

If all of $x_0$, ... , $x_6$ are zero, then $x_7=x_8=0$.

From the computation above, we see that the fixed locus of $g$ is the disjoint union of $\mathbb{P}^2\cap X$ and $\mathbb{P}^3\cap X$. A similar proof applies to the general case in the lemma below.

\begin{Lemma}
In the same setting as Theorem A, the fixed locus of a subgroup $H$ of $G$ decomposes into connected components. Each of the component is of the form $\mathbb{P}^m\cap X\subset\mathbb{P}^n$ for some $m$, where $\mathbb{P}^m$ is defined by the equations $x_{j_1}=x_{j_2}=\cdots=x_{j_{n-m}}=0$.\\
\end{Lemma}

\begin{remark}
In this paper we only consider the projective subspaces in $\mathbb{P}^n$ which are cut out by the equations $x_{j_1}=x_{j_2}=\cdots=x_{j_{n-m}}=0$.
\end{remark}

\paragraph\label{reduce} The unsimplified product on $\HT^*([X/G])$ is associative if the Bass-Quillen class associated with the sequence of schemes
$$X^{g,h}\hookrightarrow X^{gh}\hookrightarrow{X} $$
and the sequence of schemes
$$X^{g,h}\hookrightarrow X^{g}\hookrightarrow X $$
is zero~\cite{CH,H}. The Bass-Quillen class associated with a general sequence of schemes
$$Y\hookrightarrow Z\hookrightarrow S $$
is a cohomology class in $\Ext^1_{\cO_Y}(N_{Y/Z}\otimes N_{Z/S}|_Y, N_{Z/S}|_Y)$~\cite{CH, H}.

When $X$ is the Fermat hypersurface in $\mathbb{P}^{n}$, the sequence of fixed loci is of the form
$$\mathbb{P}^{l}\cap X\hookrightarrow \mathbb{P}^{m}\cap X\hookrightarrow X.$$

The following lemma shows that it suffices to study the Bass-Quillen class associated with the sequence
$$\mathbb{P}^l\hookrightarrow\mathbb{P}^m\hookrightarrow\mathbb{P}^n$$
and then restrict it to $X\cap\mathbb{P}^l$.

\begin{Lemma}\label{lemm restrict}
Let $Y\hookrightarrow S$ and $T\hookrightarrow S$ be closed embedding of smooth schemes. Assume that the intersection $W=Y\times_ST$ is smooth and transversal. Let $I_Y\subset\cO_S$ and $I_T\subset\cO_S$ be the ideal sheaves of $Y$ and $T$. Further we assume that $I_Y\cap I_T=I_Y I_T$. Then the normal bundle $N_{W/T}$ is the normal bundle $N_{Y/S}$ restricted to $W$.
\end{Lemma}

\begin{Proof}
We have a short exact sequence
$$0\rightarrow I_Y\rightarrow \cO_S\rightarrow \cO_S/I_Y\rightarrow 0$$
 on $S$.
 Tensor the exact sequence above with $\cO_T=\cO_S/I_T$. Since all the schemes are smooth and the intersection $W=Y\times_ST$ is transversal, we have $\cO_S/I_Y\otimes_{\cO_S}\cO_S/I_T$ is the structure sheaf $\cO_W$ of $W$. We obtain the following sequence

 $$I_Y|_T=I_Y\otimes\cO_S/I_T\rightarrow \cO_T\rightarrow\cO_W\rightarrow 0.$$

 Equivalently we have
 $$I_Y/I_YI_T\rightarrow\cO_S/I_T\rightarrow\cO_S/(I_Y+I_T)\rightarrow 0. $$
 The kernel of the map $\cO_S/I_T\rightarrow\cO_S/(I_Y+I_T)$ is $(I_Y+I_T)/I_T\cong I_Y/(I_Y\cap I_T)$. Due to the assumption $I_Y\cap I_T=I_YI_T$, we see that the sequence above is exact. The sequence shows that $I_Y|_T$ is the ideal sheaf of $W$ in $T$. Therefore the conormal bundle $N^{\vee}_{W/T}$ is $I_Y|_T/(I_Y|_T)^2$ which is $N^{\vee}_{Y/S}|_W$.
\end{Proof}

Set $Y$, $S$, and $T$ to be $\mathbb{P}^m$, $\mathbb{P}^n$, and $X$. The assumptions in the lemma above hold when $m\geq1$. Therefore the normal bundle of $\mathbb{P}^m\cap X\hookrightarrow X$ is the normal bundle of $\mathbb{P}^m\hookrightarrow\mathbb{P}^n$ restricted to $X\cap\mathbb{P}^m$. Similarly, one can show that the Bass-Quillen class associated with the sequence
$$\mathbb{P}^l\cap X\hookrightarrow\mathbb{P}^m\cap X\hookrightarrow\mathbb{P}^n\cap X=X$$
is the Bass-Quillen class associated with the sequence
$$\mathbb{P}^l\hookrightarrow\mathbb{P}^m\hookrightarrow\mathbb{P}^n$$
restricted to $X\cap\mathbb{P}^l$ when $l\geq1$. when $l=0$, the Bass-Quillen class vanishes on a set of points.

To study the Bass-Quillen class, we need to study the normal bundle of the map $\mathbb{P}^m\hookrightarrow\mathbb{P}^n$.
\begin{Lemma}\label{lem normal}
Let $M\subset\mathbb{P}^n$ be a complete intersection\footnote{This holds in great generality (for reducible or non-reduced $M$) if we define $N_{M/\mathbb{P}^n}$ to be $(\mathcal{I}_M/\mathcal{I}_M^2)^\vee$ although we only need the simplest case: a liner subspace of $\mathbb{P}^m$.} of irreducible polynomials $f_1$, ... , $f_j$ of degree $d_1$, ... , $d_j$. Then the normal bundle $N_{M/\mathbb{P}^n}$ is $$\bigoplus_{l=1}^{j}\cO_{M}(d_l).$$
\end{Lemma}
\begin{Proof}
One can check this by Koszul complex.
\end{Proof}

\begin{Proposition}\label{Prop bass-quillen}
The Bass-Quillen class associated to $$\mathbb{P}^l\hookrightarrow\mathbb{P}^m\hookrightarrow\mathbb{P}^n$$ is zero.
\end{Proposition}
\begin{Proof}
The normal bundle $N_{\mathbb{P}^l/\mathbb{P}^n}$ is $\cO(1)^{\oplus(n-l)}$ due to Lemma~\ref{lem normal}. The Bass-Quillen class is an element in $\Ext^1(N_{\mathbb{P}^l/\mathbb{P}^m}\otimes N_{\mathbb{P}^m/\mathbb{P}^n}|_{\mathbb{P}^l},N_{\mathbb{P}^m/\mathbb{P}^n}|_{\mathbb{P}^l})$. This Ext group is $(m-l)(n-m)(n-m)$ copies of $H^1(\mathbb{P}^l,\cO(-1))$. We know the cohomology $H^1(\mathbb{P}^l,\cO(-1))$ vanishes.
\end{Proof}

\begin{Proof}[Proof of Theorem A]
This is due to~(\ref{reduce}), Proposition~\ref{Prop bass-quillen}, and Lemma~\ref{lemm restrict}.
\end{Proof}

\section{Classification of the group elements}
In this section, we study the group $G=(\mathbb{Z}/5\mathbb{Z})^3$. The elements of $G$ are divided into four types. The fixed loci of different types of elements have different dimensions. We compute the fixed loci and the ages $\iota(g)$ for all $g\in G$. The classification will be used when we study Fermat quintic orbifold.

\paragraph Let $X$ be the Fermat quintic in $\mathbb{P}^4$ defined by $x^5+y^5+z^5+s^5+t^5=0$, where $[x:y:z:s:t]$ is the homogenous coordinates on $\mathbb{P}^4$. The group $G=(\mathbb{Z}/5\mathbb{Z})^3$ acts on $X$. An element of $G$ is of the form $g=(\zeta^{-a-b-c},\zeta^{a},\zeta^{b},\zeta^{c},\zeta^0)$, where $\zeta$ is the root of unity and $a$, $b$, $c\in\{0,1,2,3,4\}$. The group action is defined by

$$g\cdot[x:y:z:s:t]=[\zeta^{-a-b-c}x:\zeta^ay:\zeta^bz:\zeta^cs:t].$$

We classify the elements of $G$. We define the four different types of elements in $G$ as follows.

\begin{itemize}
\item Type one. There is only one element $(1,1,1,1,1)\in G$ which is of type one.

\item Type two. A nontrivial element $g=(\zeta^{-a-b-c},\zeta^{a},\zeta^{b},\zeta^{c},\zeta^0)$ is of type two if three of the following five numbers $-a-b-c$, $a$, $b$, $c$, $0$ are equal. For example $(\zeta,\zeta^4,1,1,1)$ and $(\zeta,\zeta,\zeta,\zeta^2,1)$ are of type two. There are 40 of them.

\item Type three. A nontrivial element $g=(\zeta^{-a-b-c},\zeta^{a},\zeta^{b},\zeta^{c},\zeta^0)$ is of type three if two of the following five numbers $-a-b-c$, $a$, $b$, $c$, $0$ are equal and the element is not type two. For example $(\zeta,\zeta,\zeta^4,\zeta^4,1)$ and $(1,\zeta,\zeta,\zeta^3,1)$ are of type three. There are 60 of them.

\item Type four. An element $g=(\zeta^{-a-b-c},\zeta^{a},\zeta^{b},\zeta^{c},\zeta^0)$ is of type four if all the following numbers $-a-b-c$, $a$, $b$, $c$, $0$ are different. For example $(\zeta,\zeta^2,\zeta^3,\zeta^4,1)$ is of type four. There are 24 of them.

\end{itemize}

\paragraph{\bf The fixed locus.}
The dimension of the fixed locus is completely determined by the numbers of $-a-b-c$, $a$, $b$, $c$, and $0$ that are equal to each other. The computation of the fixed locus of an element has been done in~(\ref{fixed locus}). The fixed locus $X^g$ can be classified according to the types of $g\in G$.
\begin{itemize}
\item Type one. The fixed locus is $X$.

\item Type two. The fixed locus is a genus $6$ curve in $X\cap\mathbb{P}^2$. For example, when $g=(\zeta,\zeta^4,1,1,1)$, the fixed locus is $[0:0:z:s:t]\subset X\cap\mathbb{P}^2$, where $z^5+s^5+t^5=0$. Similarly, when  $g=(\zeta,\zeta,\zeta,\zeta^2,1)$, the fixed locus is $[x:y:z:0:0]\subset X\cap\mathbb{P}^2$, where $x^5+y^5+z^5=0$.

\item Type three. The fixed locus is zero-dimensional. For example, when $g=(\zeta,\zeta,\zeta^4,\zeta^4,1)$, the fixed locus is a set of ten points $[1:-\zeta^j:0:0:0]$ and $[0:0:1,-\zeta^j:0]$.

\item Type four. No fixed locus.

\end{itemize}

\paragraph{\bf The age $\iota(g)$ of an group element $g$.}
For a general orbifold $[X/G]$, the age $\iota(g,U)$ defined in~\cite{CR} and recalled in \ref{age} is a non-negative rational number that will be used in this paper. When $X^g$ decomposes into connected components, the number $\iota(g)$ can be different on each component $U$. For simplicity, we write $\iota(g)$ instead of $\iota(g,U)$ when $U$ is clear from the context. When the orbifold is Calabi-Yau, then $\iota(g)$ is an integer~\cite{FG}. Age plays an important role in the product on the orbifold polyvector field~\cite{CH} and in the orbifold singular cohomology~\cite{CR} as well.

We can compute the age $\iota(g)$ of $g\in G$ according to its type in the case of the Fermat quintic orbifold. Type one and four are trivial. In the case of type three, one can compute the age by straightforward computations using the definition of age. For type two we may use the following identity in~\cite{FG}: $\iota(g,U)+\iota(g^{-1},U)=\mathrm{codim}(U,X)$, as the conjugation from the $S_5$ action on $X$ flips $g$ and $g^{-1}$ and preserves $U$, we see that $\iota(g,U)=\iota(g^{-1},U)=\mathrm{codim}(U,X)/2=1$. One can also use the definition to compute the age in the case of type two. To summarize:
\begin{itemize}
\item Type one. $\iota(g)=0$.
\item Type two. $\iota(g)=1$.
\item Type three. $\iota(g)=2$ or $1$. For example, choose $g=(\zeta,\zeta,\zeta^4,\zeta^4,1)$. The fixed locus is a set of ten points $[1:-\zeta^j:0:0:0]$ and $[0:0:1,-\zeta^j:0]$. The age of $g$ is $2$ on the first five points and is $1$ on the other five points.
\item Type four. No fixed locus. No $\iota(g)$.
\end{itemize}

\section{The simplified product}
C\u{a}ld\u{a}raru and Huang~\cite{CH} have a conjectural way to simplify the product we defined on $\HT^*([X/G])$. Let $\alpha_g$ be an element in ${H}^{p-c_g}(X^g,\wedge^q T_{X^g}\otimes\omega_g)$ and $\beta_h$ be an element in ${H}^{p'-c_h}(X^h,\wedge^{q'} T_{X^h}\otimes\omega_h)$. The simplified product of $\alpha_g$ and $\beta_h$ uses a cohomology class $\gamma_{g,h}$ introduced in Fantechi-G\"{o}ttsche's paper~\cite{FG}. We review the definition of the simplified product and show that the simplified product has an easier formula when the class $\gamma_{g,h}$ is trivial.  When $[X/G]$ is the Fermat quintic orbifold, we show that either the class $\gamma_{g,h}$ is trivial or the simplified product $\alpha_g\cdot\beta_h$ vanishes when $\gamma_{g,h}$ is not trivial.

We recall the definition of the simplified product and the unsimplified product in Sections 4 and 5 respectively. The rest of Sections 4 and 5 are devoted to proving Theorem B1. For most of the applications in Hochschild cohomology and homological mirror symmetry, it suffices to consider $\HT^*([X/G])$ rather than $\HT^*(X;G)$. Namely Theorem B2 is enough for most of the applications. The proof of Theorem B2 is much easier, so readers feel free to skip the proof of Theorem B1 in Sections 4 and 5.

\paragraph In~\cite{FG} Fantechi-G\"{o}ttsche introduced a cohomology class to study orbifold singular cohomology. This class $\gamma_{g,h}$ is the top Chern class of a vector bundle $R_{g,h}$ of rank $k=\iota(g)+\iota(h)-\iota(gh)-\codim(X^{g,h},X^{gh})$ on $X^{g,h}$.

\paragraph{\bf The simplified product.}\label{sim prod}
The simplified product is defined as follows
\begin{align*}
{H}^{p}(X^g,\wedge^q T_{X^g}\otimes\omega_g[-c_g])\otimes&{H}^{p'}(X^h,\wedge^{q'}T_{X^h}\otimes\omega_h[-c_h])\\
\downarrow\\
{H}^{p+p'}(X^{g,h},\wedge^qT_{X^g}|_{X^{g,h}}\otimes\omega_{g}|_{X^{g,h}}[-c_g]\otimes&\wedge^{q'}T_{X^h}|_{X^{g,h}}\otimes\omega_h|_{X^{g,h}}[-c_h])\\
={H}^{p+p'}(X^{g,h},\wedge^q T_{X^g}|_{X^{g,h}}\otimes \wedge^{q'}T_{X^h}|_{X^{g,h}}\otimes&\omega_g|_{X^{g,h}}[-c_g]\otimes\omega_h|_{X^{g,h}}[-c_h])\\
\cong{H}^{p+p'}(X^{g,h},\wedge^q T_{X^g}|_{X^{g,h}}\otimes \wedge^{q'}T_{X^h}|_{X^{g,h}}\otimes&\omega_{g,h}[-c_{g,h}]\otimes\wedge^rE[-r])\\
\downarrow\\
\bigoplus_{i=0}^{i+j=k}{H}^{p+p'-r+k}(X^{g,h},\wedge^{q-i}T_{X^g}|_{X^{g,h}}\otimes&\wedge^{q'-j}T_{X^h}|_{X^{g,h}}\otimes\omega_{g,h}[-c_{g,h}]\otimes\wedge^rE)\\
\downarrow\\
{H}^{p+p'-r+k}(X^{gh},\wedge^{q+q'-k+r}T_{X^{gh}}&\otimes\omega_{gh}[-c_{gh}]),
\end{align*}
%$${H}^{p}(X^g,\wedge^q T_{X^g}\otimes\omega_g[-c_g])\otimes{H}^{p'}(X^h,\wedge^{q'}T_{X^h}\otimes\omega_h[-c_h])$$
%$$\longrightarrow{H}^{p+p'}(X^{g,h},\wedge^q T_{X^g}|_{X^{g,h}}\otimes\omega_{g}|_{X^{g,h}}[-c_g]\otimes\wedge^{q'}T_{X^h}|_{X^{g,h}}\otimes\omega_h|_{X^{g,h}}[-c_h])
%$$
%$$
%={H}^{p+p'}(X^{g,h},\wedge^q T_{X^g}|_{X^{g,h}}\otimes \wedge^{q'}T_{X^h}|_{X^{g,h}}\otimes\omega_g|_{X^{g,h}}[-c_g]\otimes\omega_h|_{X^{g,h}}[-c_h])
%$$
%$$
%\cong{H}^{p+p'}(X^{g,h},\wedge^q T_{X^g}|_{X^{g,h}}\otimes \wedge^{q'}T_{X^h}|_{X^{g,h}}\otimes\omega_{g,h}[-c_{g,h}]\otimes\wedge^rE[-r])
%$$
%$$
%\longrightarrow \bigoplus_{i=0}^{i+j=k}{H}^{p+p'-r+k}(X^{g,h},\wedge^{q-i} T_{X^g}|_{X^{g,h}}\otimes \wedge^{q'-j}T_{X^h}|_{X^{g,h}}\otimes\omega_{g,h}[-c_{g,h}]\otimes\wedge^rE)
%$$
%$$
%\longrightarrow{H}^{p+p'-r+k}(X^{gh},\wedge^{q+q'-k+r}T_{X^{gh}}\otimes\omega_{gh}[-c_{gh}]),
%$$
where $E$ is the excess bundle of rank $r=c_g+c_h-c_{g,h}$ and $k$ is the rank of $R_{g,h}$. The first arrow is the naive restriction from $X^g$ and $X^h$ to $X^{g,h}$, so we call it by the {\em naive restriction and multiplication}. The isomorphism in the middle is due to the isomorphism~\cite{CH} below
$$\omega_{g}|_{X^{g,h}}[-c_g]\otimes\omega_h|_{X^{g,h}}[-c_h]\cong\omega_{g,h}[-c_{g,h}]\otimes\wedge^r E[-r].$$
The second arrow in the middle involving $k$ is the action of $\gamma_{g,h}$.
We call the last map the {\em extension} map because it is from $X^{g,h}$ to $X^{gh}$. More explanations can be found in~\cite{CH}.
\begin{Definition}
We say the class $\gamma_{g,h}$ is trivial if the rank $k=\iota(g)+\iota(h)-\iota(gh)-\codim(X^{g,h},X^{gh})$ is zero or the rank $k$ is strictly greater than the dimension of $X^{g,h}$. The class $\gamma_{g,h}$ is $1$ in the first case and $0$ in the second case.\\
\end{Definition}

\begin{remark}
The rank $k$ of $R_{g,h}$ may be greater than the dimension of $X^{g,h}$. For example, choose $g=h=(\xi^4,\xi^4,\xi^2,1,1)$ when $[X/G]$ is the Fermat quintic orbifold. Then $\iota(g)=\iota(h)=\iota(gh)=2$ and $c_g=c_h=c_{gh}=c_{g,h}=3$. The vector bundle $R_{g,h}$ has rank $2$ on $X^{g,h}$ which is a set of points.
\end{remark}

We study the simplified product when $\gamma_{g,h}$ is trivial. We need the proposition below.

\begin{Proposition}\label{Prop ineq}
When $k>\dim X^{g,h}$, $r>\dim X^{gh}$.
\end{Proposition}

\begin{Proof}
Let $d_{g}$, $d_{g,h}$ be the dimension of $X^{g}$ and $X^{g,h}$. The rank $k$ is by definition $\iota(g)+\iota(h)-\iota(gh)-d_{gh}+d_{g,h}.$
When $k>d_{g,h}$, we have $\iota(g)+\iota(h)>d_{gh}-\iota(gh)$. There is an equality
$$\iota(g)+\iota(g^{-1})=c_g$$
in~\cite{FG}. Therefore $r=c_g+c_h-c_{g,h}=\iota(g)+\iota(g^{-1})+\iota(h)+\iota(h^{-1})-c_{g,h}$.
Using the inequality above, we have
\begin{equation}\label{eq1}
r>\iota(g^{-1})+\iota(h^{-1})+d_{gh}-\iota(gh)-c_{g,h}.
\end{equation}
On the other hand, the rank
$$r=\iota(g)+\iota(h)-\iota(gh)-d_{gh}+d_{g,h}=\iota(g)+\iota(h)-\iota(gh)+c_{gh}-c_{g,h}\geq0$$
for all $g,h\in G$, i.e.,
$$\iota(g)+\iota(h)-\iota(gh)+c_{gh}\geq c_{g,h}$$
for all $g,h\in G$.
Note that $c_{gh}-\iota(gh)=\iota((gh)^{-1})$ which implies the following
$$\iota(g)+\iota(h)-\iota((gh)^{-1})\geq c_{g,h}.$$
Similarly, we have
\begin{equation}\label{eq2}
\iota(h^{-1})+\iota(g^{-1})-\iota(gh)\geq c_{h^{-1},g^{-1}}=c_{h,g}.
\end{equation}
Combining the inequalities (\ref{eq1}) and (\ref{eq2}) above, we get the desired inequality $r>d_{gh}$.
\end{Proof}

\begin{Proposition}\label{Prop trivial}
When $\gamma_{g,h}$ is trivial, the simplified product is equal to the composite map below
\begin{align*}
{H}^{p}(X^g,\wedge^q T_{X^g}\otimes\omega_g[-c_g])\otimes & {H}^{p'}(X^h,\wedge^{q'}T_{X^h}\otimes\omega_h[-c_h])\\
\downarrow\\
{H}^{p+p'}(X^{g,h},\wedge^qT_{X^g}|_{X^{g,h}}\otimes\omega_{g}|_{X^{g,h}}[-c_g]\otimes & \wedge^{q'}T_{X^h}|_{X^{g,h}}\otimes\omega_h|_{X^{g,h}}[-c_h])\\
={H}^{p+p'}(X^{g,h},\wedge^q T_{X^g}|_{X^{g,h}}\otimes \wedge^{q'}T_{X^h}|_{X^{g,h}}\otimes & \omega_g|_{X^{g,h}}[-c_g]\otimes\omega_h|_{X^{g,h}}[-c_h])\\
\cong{H}^{p+p'-r}(X^{g,h},\wedge^q T_{X^g}|_{X^{g,h}}\otimes \wedge^{q'}T_{X^h}|_{X^{g,h}}\otimes & \omega_{g,h}[-c_{g,h}]\otimes\wedge^rE)\\
\downarrow\\
{H}^{p+p'-r}(X^{gh},\wedge^{q+q'+r}T_{X^{gh}}\otimes & \omega_{gh}[-c_{gh}]).
\end{align*}
Namely, we only do the naive restriction and multiplication and then extend.
\end{Proposition}

\begin{Proof}
When the rank $k$ of $R_{g,h}$ is zero, the class $\gamma_{g,h}$ is equal to $1$. It is clear that $\gamma_{g,h}$ acts as the identity map.

When the rank $k$ of $R_{g,h}$ is greater than the dimension of $X^{g,h}$, the class is equal to $0$ rather than $1$. It is clear that $\gamma_{g,h}$ acts as the zero map. It suffices to show that the composite map in Proposition~\ref{Prop trivial} is also zero in this case. Due to Proposition~\ref{Prop ineq}, we have $r>\dim X^{gh}$ which implies the last term ${H}^{p+p'-r}(X^{gh},\wedge^{q+q'+r}T_{X^{gh}}\otimes\omega_{gh}[-c_{gh}])$ of the composite map vanishes.
\end{Proof}

When $\gamma_{g,h}$ is trivial, the product can be summarized as restriction and multiplication which is the first arrow, and then extension which is the last arrow. There is no nontrivial construction in the middle in this special case.

\paragraph{\bf The simplified product for Fermat quintic.} From now on we study the Fermat quintic case.

\begin{Proposition}~\label{Prop quintic}
In the same setting as Theorem B1, either the class $\gamma_{g,h}$ is trivial or $g$ is of type two and $h=g^{j}$, where $j=1,2,3$.
\end{Proposition}
\begin{Proof}
The class depends on $g$, $h$, and $(gh)^{-1}$. In~\cite{FG} it is shown that $\gamma_{g,h}$ has the following property $(*)$

$$\gamma_{g,h}=\gamma_{h,g}=\gamma_{g,(gh)^{-1}}=\gamma_{h,(gh)^{-1}}.$$

We study the class $\gamma_{g,h}$ case by case according to the dimension of $X^{g,h}$.

The group $G$ can be viewed as a vector space $V=G$ over $\mathbb{Z}/5\mathbb{Z}$ and $g,h$ can be viewed as vectors. If the two elements $g,h$ are linearly independent, i.e., $h$ is not in the cyclic group $<g>$ generated by $g$, then $X^{g,h}$ is zero-dimensional. If the vectors $g,h$ generate a one-dimensional subspace of $V$, i.e., $h$ is an element of the cyclic group $<g>$ generated by $g$, then $X^{g,h}$ could be a genus $6$ curve.

When $X^{g,h}$ is zero dimensional, the class $\gamma_{g,h}$ is trivial.

Consider the case when $X^{g,h}$ is a curve. Then $g$ must be an element of type two and $h$ must be of the form $g^j$, where $j=0,1,2,3,4$. Due to the property $(*)$ of the class $\gamma_{g,h}$, it suffices to consider the cases when $h=1$ and $h=g$.

Recall that the class $\gamma_{g,h}$ is the top Chern class of a vector bundle on $X^{g,h}$ of rank
$$\iota(g)+\iota(h)-\iota(gh)-\mathrm{codim}(X^{g,h},X^{gh}).$$
When $h=1$, it is clear that the rank is zero and therefore $\gamma_{g,h}$ is trivial.

When $h=g$, it is clear that $\iota(g)=\iota(g^2)=1$ because the fixed locus $X^{g}$ is a curve, and therefore $g$ is of type two. The rank is $1$.

Similarly, we see that the rank $r$ is zero when $h=g^4$ and the rank is $1$ when $h=g^j$, where $j=2$ and $3$.
\end{Proof}

\begin{Proposition}
We are in the same setting as Theorem B1. Let $\alpha_g$ be an element in ${H}^{p-c_g}(X^g,\wedge^q T_{X^g}\otimes\omega_g)$ and $\beta_h$ be an element in ${H}^{p'-c_h}(X^h,\wedge^{q'} T_{X^h}\otimes\omega_h)$. The simplified product of $\alpha_g$ and $\beta_h$ is equal to the composite map in Proposition~\ref{Prop trivial} for all $g,h\in G$.
\end{Proposition}

\begin{Proof}
Due to Propositions~\ref{Prop trivial} and~\ref{Prop quintic}, it suffices to consider the case where $g$ is of type two and $h=g^{j}$ for $j=1,2,3$. It suffices to show that both the simplified product in~(\ref{sim prod}) and the composite map in Proposition~\ref{Prop trivial} are zero in this case. Note that $X^{g,h}=X^{gh}$ is a genus $6$ curve and the rank $r=c_g+c_h-c_{g,h}=2$ in this case.
A map of the form
$$
{H}^{*}(X^{g,h},\wedge^* T_{X^g}|_{X^{g,h}}\otimes \wedge^{*}T_{X^h}|_{X^{g,h}}\otimes\omega_{g,h}\otimes\wedge^2 E)
$$
$$
\longrightarrow{H}^{*}(X^{gh},\wedge^{*+2}T_{X^{gh}}\otimes\omega_{gh}).
$$
appears in the last arrow of the simplified product in~(\ref{sim prod}) and in the last arrow of the composite map in Proposition~\ref{Prop trivial}. The vector bundle $\wedge^{*+2}T_{X^{gh}}$ is zero because $2$ is greater than the dimension of $X^{gh}$.
\end{Proof}

\section{Proof of Theorem B1}
We review the definition of the unsimplified product and prove Theorem B1.

\paragraph{\bf The unsimplified product.}
The unsimplified product is defined in a similar way~\cite{CH}. We do the derived restriction and multiplication first and then extend
$$
	\mathrm{D}(\mathbb{L}_{\widetilde{X^g}}/ X)\otimes\mathrm{D}(\mathbb{L}_{\widetilde{X^h}}/ X) \rightarrow \mathrm{D}(\mathbb{L}_{\widetilde{X^g}}\times^R_{X}\mathbb{L}_{\widetilde{X^h}}/ X)\stackrel{\cong}{\rightarrow} \mathrm{D}(\mathbb{L}_{(\widetilde{X^g}\times_{X} \widetilde{X^h})}/ X)
$$
$$
 \stackrel{{\mathbb{L}_{m}}_*}{\longrightarrow} \mathrm{D}(\mathbb{L}_{\widetilde{X^{gh}}}/ X),
$$
where the maps above are explained in~\cite{CH}. The only difference is that the first arrow in the simplified product is the naive restriction and the first arrow in the unsimplified product is the derived restriction which could have more terms
$${H}^{p-c_g}(X^g,\wedge^q T_{X^g}\otimes\omega_g)\otimes{H}^{p'-c_h}(X^h,\wedge^{q'}T_{X^h}\otimes\omega_h)$$
$$
\longrightarrow\bigoplus_{i=0}^{r}{H}^{p+p'-i-c_{g,h}}(X^{g,h},\wedge^q T_{X^g}|_{X^{g,h}}\otimes \wedge^{q'}T_{X^h}|_{X^{g,h}}\otimes\omega_{g,h}\otimes\wedge^iE).
$$
Then we use the same extension map as before, so the output of the unsimplified product lands in
$$\longrightarrow \bigoplus_{i=0}^{r}H^{p+p'-c_{gh}-i}(X^{gh},\wedge^{q+q'+i}T_{X^{gh}}\otimes\omega_{gh}).$$

\begin{Proof}[Proof of Theorem B1.]
Let $\alpha_g$ be an element in ${H}^{p-c_g}(X^g,\wedge^q T_{X^g}\otimes\omega_g)$ and $\beta_h$ be an element in ${H}^{p'-c_h}(X^h,\wedge^{q'} T_{X^h}\otimes\omega_h)$. To prove the two products agree, it suffices to show that
$$\alpha_g\cdot\beta_h\in H^{p+p'-c_{gh}-r}(X^{gh},\wedge^{q+q'+r}T_{X^{gh}}\otimes\omega_{gh})$$
for the unsimplified product. Namely, the unsimplified product of $\alpha_g$ and $\beta_h$ lands only in one of the direct summand, where $i$ can only be $r$, of the big direct sum
$$\bigoplus_{i=0}^{r}H^{p+p'-c_{gh}-i}(X^{gh},\wedge^{q+q'+i}T_{X^{gh}}\otimes\omega_{gh}).$$
We prove the statement above according to the dimension of the fixed locus $X^{g,h}$.

When $X^{g,h}$ is zero dimensional, we know that $p+p'-c_{g,h}-i$ must be zero. We also have $p-c_g\geq0$ and $p'-c_h\geq0$. We conclude $i+c_{g,h}=p+p'\geq c_g+c_h$. However, $i$ is an integer from $0$ to $r=c_g+c_h-c_{g,h}$ which forces $i$ to be $r$ in this case.

When $X^{g,h}$ is not zero dimensional, $h$ must be of the form $g^j$, where $j=0,1,2,3,4$. The fixed locus is a genus $6$ curve $C$.

When $h=g^0=(1,1,1,1,1)$ is the identity, it is easy to see that the rank $r=c_g+c_h-c_{g,h}$ is zero.

When $h=g^{-1}$, $c_g=c_h=c_{g,h}=2$ and $c_{gh}=0$. We look at the unsimplified product

$$H^{p-2}(C,\wedge^qT_{C}\otimes\omega_C)\otimes H^{p'-2}(C,\wedge^{q'}T_{C}\otimes\omega_C) $$
$$\rightarrow\bigoplus_{i=0}^{2} H^{p+p'-2-i}(C,\wedge^qT_{C}\otimes\wedge^{q'}T_{C}\otimes\wedge^i E\otimes\omega_{C}) $$
$$\rightarrow \bigoplus_{i=0}^{2}H^{p+p'-i}(X,\wedge^{q+q'+i}T_{X}),$$
where $E$ is the normal bundle $N_{C/X}$.

Recall that we need to prove that all the maps vanish when $i\neq2$. We know that $1\geq p-2\geq0$ and $1\geq p'-2\geq0$ from the first line in the product and $1\geq p+p'-2-i\geq 0$ from the second line in the product. We conclude that $i$ can not be $0$ immediately.

Consider the case when $i=1$. Due to the same inequality above, we can conclude that $p$ and $p'$ must be $2$ in this case. The dimensions of $\HT^*(X)$ is well-known and nonzero terms in $\HT^{*,3}(X)$ are $\HT^{0,3}(X)$ and $\HT^{3,3}(X)$. This implies $q$ and $q'$ must be $1$. The last map in the product above is induced by a map of vector bundles $\wedge^{q}T_{C}\otimes\wedge^{q'}T_C\otimes\wedge^i E\rightarrow \wedge^{q}T_X|_C\otimes\wedge^{q'}T_X|_C\otimes\wedge^{i}T_X|_C\rightarrow\wedge^{q+q'+i} T_X|_C$~\cite{CH}. When $q=q'=1$, the map of vector bundles must vanish because $T_C$ is rank $1$ and $\wedge^2T_C=0$. We complete the proof that $i$ must be $2$ when $h=g^{-1}$.

When $h=g^j$ for $j=1,2,3$, we can assume $j=1$ and $h=g$ without losing generality. One can check that the proof for the case when $j=2,3$ is similar to the proof below. Under this assumption, $c_g=c_h=c_{g,h}=c_{gh}=2$. The rank $r$ is also equal to $2$. We look at the unsimplified product
$$H^{p-2}(C,\wedge^qT_{C}\otimes\omega_C)\otimes H^{p'-2}(C,\wedge^{q'}T_{C}\otimes\omega_C) $$
$$\rightarrow\bigoplus_{i=0}^{2} H^{p+p'-2-i}(C,\wedge^qT_{C}\otimes\wedge^{q'}T_{C}\otimes\wedge^i E\otimes\omega_{C}) $$
$$\rightarrow \bigoplus_{i=0}^{2}H^{p+p'-2-i}(C,\wedge^{q+q'+i}T_{C}\otimes\omega_C),$$
where $E=N_{C/X}$ is the normal bundle in this case.

Recall that we need to prove that all the maps vanish when $i\neq2$. We know that $1\geq p-2\geq0$ and $1\geq p'-2\geq0$ from the first line in the product and $1\geq p+p'-2-i\geq 0$ from the second line in the product. We conclude that $i$ can not be $0$ immediately.

Consider the case when $i=1$. Due to the same inequality above, we can conclude that $p$ and $p'$ must be $2$ in this case.

The last line above shows that $1\geq q+q'+i=q+q'+1\geq0$ which shows that $q$ and $q'$ must be zero in this case.

In this case, the map we are looking at is
$$H^{0}(C,\omega_C)\otimes H^{0}(C,\omega_C) $$
$$ \rightarrow H^{1}(C, E\otimes\omega_{C}) $$
$$\rightarrow H^{1}(C,T_{C}\otimes\omega_C)=H^1(C,O_C),$$
The second arrow above is induced by a map $E\rightarrow T_{X^{gh}}|_{X^{g,h}}=T_C$~\cite{CH}. We can show that the map is zero as follows. Denote $T_X|_C$ by $V$. Then $E=\frac{V}{V^g+V^h}=\frac{V}{V^g}$ in this case. The tangent bundle $T_{X^{gh}}$ is naturally considered as a quotient space $V_{gh}$ of $V$~\cite{ACH}, not a subspace of $V$. The map
$$E=\frac{V}{V^g+V^h}\rightarrow V\rightarrow V_{gh}  $$
is defined by the formula~\cite{CH}
$$v\rightarrow v-g\cdot v.$$
Note that $V_{gh}=\frac{V}{<v-gh\cdot v>}$. In the case when $h=g$, we have $v-gh\cdot v=v-g^2\cdot v$.

The relation above means $v=gh\cdot v=g^2\cdot v$ for a vector $v\in V_{gh}$. Then
$$v\rightarrow v-g\cdot v=v-g(g^2\cdot v)=v-g(g^2(g^2\cdot v))=v-v=0$$
which shows that the map $E\rightarrow V_{gh}$ is zero. We can conclude that the map
$$H^{1}(C, E\otimes\omega_{C})\rightarrow H^{1}(C,T_{X^{g,h}}\otimes\omega_C)=H^1(C,T_C\otimes\omega_C)$$
is zero in this case.

As a consequence, the only possible non-vanishing product lands in the direct summand where $i=2$.

Note that it is not hard to show that the product is also zero when $i=2$ and $j=1,2,3$ by direct computation.
\end{Proof}

\section{The multiplicative bigrading}
Before we prove Theorem B2, we need to study the multiplicative bigrading and prove Theorem C. We do not need $[X/G]$ to be the Fermat hypersurface in this section. We assume $[X/G]$ is an arbitrary Calabi-Yau orbifold of dimension $n$.
\begin{Proof}[Proof of Theorem C.]
The dualizing sheaf $\omega_g$ of the map $X^g\hookrightarrow X$ is the top exterior power of the normal bundle $N_{X^g/X}$~\cite{Hart}. There is a short exact sequence of sheaves on $X^g$
$$0\rightarrow T_{X^g}\rightarrow T_X\rightarrow N_{X^g/X}\rightarrow0.$$
It implies that $\omega_{X^g}^{\vee}\otimes\omega_g\cong \omega^{\vee}_X|_{X^g}$ by taking the top exterior power of the sequence. Since $[X/G]$ is Calabi-Yau, the canonical bundle $\omega_X$ is trivial. Therefore $\omega_g\cong\omega_{X^g}$.
Due to the nondegenerate pairing
$$\wedge^{q}\Omega_{X^{g}}\otimes\wedge^{d_g-q}\Omega_{X^g}\rightarrow\omega_{X^g}\cong\omega_g,$$
we can identify $\wedge^qT_{X^g}\otimes\omega_g$ with $\wedge^{d_g-q}\Omega_{X^g}$, where $d_g$ is the dimension of $X^g$. Therefore
$$H^{p-\iota(g)}(X^g,\wedge^{q-c_g+\iota(g)}T_{X^g}\otimes\omega_g) $$
is isomorphic to
$$H^{p-\iota(g)}(X^g,\wedge^{d_g-q+c_g-\iota(g)}\Omega_{X^g}).$$
Note that $d_g+c_g=n$, so $$H^{p-\iota(g)}(X^g,\wedge^{d_g-q+c_g-\iota(g)}\Omega_{X^g})\cong H^{p-\iota(g)}(X^g,\wedge^{n-q-\iota(g)}\Omega_{X^g}).$$
Taking the sum over $g\in G$, we obtain the first result $\HT^{q,p}(X;G)\cong\HOmega^{n-q,p}(X;G)$. Moreover, all the isomorphisms are functorial, hence equivariant under the $G$ action on $X$, which implies that the isomorphism $\HT^{q,p}(X;G)\cong\HOmega^{n-q,p}(X;G)$ is compatible with the $G$ action, taking invariants we get the second result $\HT^{q,p}([X/G])\cong\HOmega^{n-q,p}([X/G])$.
\end{Proof}

\section{Proof of Theorem B2 and Theorem D}
In this section, we assume that $[X/G]$ is a Calabi-Yau Fermat hypersurface orbifold. We put the dimensions of $\HT^{p,q}([X/G])$ in the form of a diamond. We study this diamond and compute it explicitly as an example when $[X/G]$ is the Fermat quintic orbifold. Then we prove Theorem B2 and Theorem D.

%We compute the product on the orbifold polyvector field. We explain and show that it matches the product of the state space of the FJRW theory of the Fermat polynomial. This is expected due to mirror symmetry.

\paragraph\label{7.1} In mirror symmetry, the mirror of the Calabi-Yau Fermat hypersurface orbifold $[X/G]$ is $X$. Therefore we expect a natural identification 
\[
\HT^{q,p}([X/G])\cong\HH^{q,p}([X/G])\cong\HOmega^{q,p}(X).
\]
as bigraded vector spaces.
This identification is due to Theorem C and the homological mirror symmetry conjecture. Below we explain a concrete proof of the identification which can be found in the literature.

First, we apply Orlov's theorem of derived equivalence of categories. Orlov's theorem~\cite{Hir} says that there is a canonical equivalence of categories 
$$D^b(X)\cong MF_{gr}(\mathbb{A}^{n+1},\mathbb{Z}/d\mathbb{Z},\Sigma_{i=0}^n\,x_i^d),$$
where the left hand side is the derived category of $X$ and the right hand side is the $\mathbb{Z}$-graded matrix factorization category with the $\mathbb{Z}/d\mathbb{Z}$-action on the coordinates $x_i$. The Hochschild cohomology is a categorical invariant, so we obtain the isomorphism below
$$\HH^{*}(X)\cong \HH^*(D^b(X))\cong \HH^{*}(MF_{gr}(\mathbb{A}^{n+1},\mathbb{Z}/d\mathbb{Z},\Sigma_{i=0}^n \,x_i^d)).
$$

Note the fact that the derived category and matrix factorization category of a $G$-quotient can be viewed as the $G$-invariant of the original category ~\cite{gr,mckay}, so this isomorphism can be upgraded to $G$-equivariant version:

\begin{align*}
\HH^{*}([X/G])&\cong \HH^*(D^b([X/G]))\cong \HH^*(D^b(X)^G)\\
&\cong \HH^{*}(MF_{gr}(\mathbb{A}^{n+1},\mathbb{Z}/d\mathbb{Z},\Sigma_{i=0}^n \,x_i^d)^G)\\
&\cong\HH^{*}(MF_{gr}(\mathbb{A}^{n+1},G\times \mathbb{Z}/d\mathbb{Z},\Sigma_{i=0}^n \,x_i^d))    \\
\end{align*}

Then it has been proven that the following two are identified as bigraded vector spaces~\cite{kr}

$$
\HH^{q,p}(MF_{gr}(\mathbb{A}^{n+1},G\times \mathbb{Z}/d\mathbb{Z},\Sigma_{i=0}^n \,x_i^d))\cong H^{q,p}_{FJRW}(W),
$$
where the right hand side is the state space of the FJRW theory of the Fermat polynomial $W=\sum_{i=0}^{n}x_i^d$.

It has been proven~\cite{Ch11} that the state space of the FJRW theory is identified with $\HOmega^{q,p}(X)$ as bigraded vector spaces. Putting all the identifications above together, we get the desired result. In fact, taking direct sum over all $(q,p)$, homological mirror symmetry predicts that the identifications above should be isomorphisms of bigraded algebras, where the product on $H^*_{FJRW}(W)$ and $H^*(X,\mathbb{C})=\displaystyle{\bigoplus_{q+p=*}\HOmega^{q,p}(X)}$ are the quantum products.

\paragraph\label{7.2} In paragraphs~\ref{7.2} and ~\ref{diamond}, we compute the dimensions $\HT^{q,p}([X/G])$ explicitly as an example to illustrate the identification above.
%As mentioned in the introduction, the bidegree below on $\HT^*(X;G)$
%\[
% \HT^{q,p}(X;G)=\bigoplus_{g\in
%    G}H^{p-\iota(g)}(X^g,\wedge^{q-c_g+\iota(g)}
%  T_{X^g}\otimes\omega_g)\!\!\!\!\!,
%\]
%is preserved by the simplified product. First we compute the dimensions of the items in the big direct sum above.

Denote $H^{p-\iota(g)}(X^g,\wedge^{q-\iota(g)}\Omega_{X^g})$ by $\HOmega^{q,p}(X;g)$ and its dimension by $h^{q,p}(X;g)$. Similarly denote $H^{p-\iota(g)}(X^g,\wedge^{q-c_g+\iota(g)}T_{X^g}\otimes\omega_g)$ by $\HT^{q,p}(X;g)$ and its dimension by $\check{h}^{q,p}(X;g)$. We compute the numbers $h^{q,p}(X;g)$ and $\check{h}^{q,p}(X;g)$ for Fermat quintic orbifold. The age and the fixed locus from Sections 2 and 3 give the numbers immediately. We put the numbers into a diamond. See the pictures below.
\[
\begin{tikzcd}[row sep=tiny, column sep=tiny]
  &   &     & 1 &     &   &   & & & &   &   &   & & & &   & & & & & \\
  &   & 0   &   & 0   &   &   & & & &   &   &   & & & &   & & & & & \\
  & 0 &     & 1 &     & 0 &   & & & &   & 1 &   & & & & 1\times 5=5 & & & & & \\
1 &   & 101 &   & 101 &   & 1 & & & & 6 &   & 6 & & & &   & & & & & \\
  & 0 &     & 1 &     & 0 &   & & & &   & 1 &   & & & & 1\times 5=5 & & & & & \\
  &   & 0   &   & 0   &   &   & & & &   &   &   & & & &   & & & & & \\
  &   &     & 1 &     &   &   & & & &   &   &   & & & &   & & & & &
\end{tikzcd}
\]
The three diamonds above from the left to the right correspond to the numbers $h^{q,p}(X;g)$ when $g$ is of type one, two, and three respectively. There are $60$ type three elements whose fixed loci are points denoted by $\{*\}$. For each of the type three elements, the fixed locus is a set of ten points. The element $g$ acts transitively on five of them and transitively on the other five of them. The age of $g$ is $1$ or $2$ depending on the points, so the cohomology $H^0(\{*\},\mathbb{C})$ of the fixed locus has degree $(1,1)$ and $(2,2)$ respectively.

\[
\begin{tikzcd}[row sep=tiny, column sep=tiny]
  &   &     & 1   &     &   &   & & & &   &   &   & & &   &   &   & & & & \\
  &   & 0   &     & 0   &   &   & & & &   &   &   & & &   &   &   & & & & \\
  & 0 &     & 101 &     & 0 &   & & & &   & 6 &   & & &   &   &   & & & & \\
1 &   & 1   &     & 1   &   & 1 & & & & 1 &   & 1 & & & 1\times 5 &   & 1\times 5 & & & & \\
  & 0 &     & 101 &     & 0 &   & & & &   & 6 &   & & &   &   &   & & & & \\
  &   & 0   &     & 0   &   &   & & & &   &   &   & & &   &   &   & & & & \\
  &   &     & 1   &     &   &   & & & &   &   &   & & &   &   &   & & & &
\end{tikzcd}
\]
The three diamonds above from the left to the right correspond to the numbers $\check{h}^{q,p}(X;g)$ when $g$ is of type one, two, and three respectively. There are $60$ type three elements whose fixed loci are points denoted by $\{*\}$. The cohomology $H^0(\{*\},\mathbb{C})$ of the fixed locus has degree $(1,2)$ and $(2,1)$ respectively due to the same reason above.

\paragraph\label{diamond} Let $h\in G$ be an element of the group $G$ and $\alpha_g\in\HT^{q,p}(X;g)$ be a class indexed by $g$. Then $h\cdot\alpha_{g}\in\HT^{q,p}(X;hgh^{-1})$ is a class indexed by $hgh^{-1}$. When the group $G$ is abelian, the group $G$ acts on each direct summand $\HT^{q,p}(X;g)$ of $\HT^{q,p}(X;G)$ individually. Denote the dimensions of $\HT^{q,p}(X;g)^G$
and $\HOmega^{q,p}(X;g)^G$ by $\check{h}^{q,p}(X;g)^G$ and $h^{q,p}(X;g)^G$ respectively. We compute the dimensions in the following.

In the proof of Theorem C, the nondegenerate pairing
$$\wedge^{q}\Omega_{X^{g}}\otimes\wedge^{d_g-q}\Omega_{X^g}\rightarrow\omega_{X^g}\cong\omega_g$$
identifies
$$\HT^{q,p}(X;g)=H^{p-\iota(g)}(X^g,\wedge^{q-c_g+\iota(g)}T_{X^g}\otimes\omega_g) $$
with
$$\HOmega^{n-q,p}(X;g)=H^{p-\iota(g)}(X^g,\wedge^{d_g-q+c_g-\iota(g)}\Omega_{X^g}),$$
where $d_g$ is the dimension of $X^g$. Their dimensions are related by the equality $h^{n-q,p}(X;g)=\check{h}^{q,p}(X;g)$.

In the case of the Fermat quintic orbifold, the pairing above is compatible with the group action. Therefore, the corresponding $G$-invariants $\HT^{q,p}(X;g)^G$
and $\HOmega^{3-q,p}(X;g)^G$ are naturally isomorphic. To compute $\check{h}^{q,p}(X;g)^G$, it suffices to compute $h^{3-q,p}(X;g)^G$. We compute the numbers when $g$ is of type four, three, two, and one respectively.
\begin{itemize}
\item There are 24 elements of type four of the group $G$. The fixed locus is empty in this case.

\item There are 60 elements of type three of the group $G$. As an example we choose $g=(\zeta,\zeta,\zeta^4,\zeta^4,1)$. The computations for general type three elements are similar. The fixed locus is a set of ten points $[1:-\zeta^j:0:0:0]$ and $[0:0:1,-\zeta^j:0]$. The age $\iota(g)$ is $2$ on the first five points and $1$ on the other five points. Therefore the dimension

$$\check{h}^{q,p}(X;g)=\dim\HT^{q,p}(X;g)=\dim H^{p-\iota(g)}(X^g,\wedge^{q-c_g+\iota(g)}T_{X^g}\otimes\omega_g)$$

is five when $(q,p)=(1,2)$ or $(2,1)$ and zero otherwise.

The group $G$ acts transitively on the first five points $[1:-\zeta^j:0:0:0]$ and transitively on the other five points $[0:0:1,-\zeta^j:0]$ respectively. The $G$-invariant dimension $\check{h}^{q,p}(X;g)^G$ is one when $(q,p)=(1,2)$ or $(2,1)$ and zero otherwise.

\item There are 40 elements of type two of the group $G$. The age $\iota(g)$ is one. The fixed locus is a genus $6$ curve $C\subset X\cap\mathbb{P}^2$. It has been shown above that $h^{1,1}(X;g)=\check{h}^{2,1}(X;g)=1$, $h^{2,2}(X;g)=\check{h}^{1,2}(X;g)=1$, $h^{2,1}(X;g)=\check{h}^{1,1}(X;g)=6$, $h^{1,1}(X;g)=\check{h}^{2,2}(X;g)=6$ and zero otherwise. It is easy to see that the group $G$ acts trivially on $H^{1-1}(X^g,\wedge^{1-1}\Omega_C)=H^0(C,\cO_C)$. This implies that the $G$-invariant dimension $h^{1,1}(X;g)^G$ is one and similarly $h^{2,2}(X;g)^G$ is one due to Serre duality. In the next paragraph, we show that the $G$ action on $H^{1-1}(X,\wedge^{2-1}\Omega_C)=H^0(X,\Omega_C)$ has no invariant. This implies that the $G$-invariant dimension $h^{2,1}(X;g)^G$ is zero and $h^{1,2}(X;g)^G$ is zero due to Serre duality.

We show that $H^0(C,\Omega_C)$ has no invariant under the $G$ action. We take $g=(\zeta,\zeta,\zeta,\zeta^2,1)$ and the fixed locus $C$ is $[x:y:z:0:0]\subset X\cap\mathbb{P}^2$, where $x^5+y^5+z^5=0$. The computation for general type three elements is similar. A basis of differential forms of this curve is given by the formula below~\cite{Hi}. Let $y_2=\frac{y}{x}$ and $y_3=\frac{z}{x}$. Then

$$\theta_{r,\alpha}=\frac{y_2^rdy_2}{y_3^\alpha},$$

where $0\leq\alpha\leq4$ and $0\leq r\leq\alpha-2$, form a basis of $H^0(C,\Omega_C)$. One can check directly that none of them is invariant under the group action.

\item There is one element of type one. The fixed locus is the entire $X$. The Hodge diamond is shown above. We want to show that the Hodge diamond after taking the $G$-invariant is of the following form
\[
\begin{tikzcd}[row sep=tiny, column sep=tiny]
  &   &     & 1   &     &   &   & \\
  &   & 0   &     & 0   &   &   & \\
  & 0 &     & 1 &     & 0 &   & & \\
1 &   & 1   &     & 1   &   & 1 & \\
  & 0 &     & 1 &     & 0 &   & & \\
  &   & 0   &     & 0   &   &   & \\
  &   &     & 1.   &     &   &   &
\end{tikzcd}
\]

To show this, we apply Orlov's theorem of derived equivalence of categories~\cite{Hir} 
$$D^b(X)\cong MF_{gr}(\mathbb{A}^5,\mathbb{Z}/5\mathbb{Z},\Sigma_{i=0}^4\,x_i^5),$$
where the left hand side is the derived category of $X$ and the right hand side is the $\mathbb{Z}$-graded matrix factorization category with the $\mathbb{Z}/5\mathbb{Z}$-action on the coordinates $x_i$. The Hochschild cohomology is a categorical invariant, so we obtain the isomorphism below
$$\HH^*(X)\cong \HH^*(MF_{gr}(\mathbb{A}^5,\mathbb{Z}/5\mathbb{Z},\Sigma_{i=0}^4 \,x_i^5)).$$

It is known ~\cite{cal,twisted} that the Hochschild cohomology of the equivariant matrix factorization $MF(\mathbb{A}^{5},\mathbb{Z}/5\mathbb{Z},\sum_{i=0}^{4}x_i^5)$ is isomorphic to a sum as follows (we refer the readers to \cite{twisted} for the product structure and this will not be used in our paper)

$$\HH^*(X)\cong(\bigoplus_{g\in\mathbb{Z}/5\mathbb{Z}}\mathsf{Jac}(Y^g,W|_{Y^g})\otimes\omega_g)^{\mathbb{Z}/5\mathbb{Z}},$$
where $Y$ is the affine space $\mathbb{A}^5$, $W=\sum_{i=0}^{4}x_i^5$ and $\omega_g$is the dualizing sheaf of the embedding $Y^g\hookrightarrow Y$. Note that $Y^g$ is a point when $g\in\mathbb{Z}/5\mathbb{Z}$ is nontrivial. We know that the dimensions of the Hochschild cohomology\footnote{Note that the Hodge diamond for Hochschild homology and cohomology are related by a flip which corresponds to contraction with the holomorphic volume form.} of $X$ are of the form
\[
\begin{tikzcd}[row sep=tiny, column sep=tiny]
  &   &     & 1   &     &   &   & \\
  &   & 0   &     & 0   &   &   & \\
  & 0 &     & 101 &     & 0 &   & & \\
1 &   & 1   &     & 1   &   & 1 & \\
  & 0 &     & 101 &     & 0 &   & & \\
  &   & 0   &     & 0   &   &   & \\
  &   &     & 1.  &     &   &   &
\end{tikzcd}
\]
Under the isomorphism above, the odd degree part of the Hochschild cohomology $\HH^*(X)$ corresponds to the twisted Jacobi ring. $\mathsf{Jac}(Y^g,W|_{Y^g})\otimes\omega_g$ indexed by the four nontrivial elements of $\mathbb{Z}/5\mathbb{Z}$. The even degree part of $\HH^*(X)$ is the vertical line in the diamond above. It corresponds to the Jacobi ring $\mathsf{Jac}(\mathbb{A}^5,W)$ of the Fermat polynomial $W$. The numbers $1$, $101$, $101$, and $1$ in the diamond of the Hochschild cohomology correspond to the numbers of monomials of degree $0$, $5$, $10$, and $15$ in the Jacobi ring.

We add the group $G=(\mathbb{Z}/5\mathbb{Z})^3$ action on both sides of the isomorphism above. One concludes that the odd degree part is invariant under the group action. In the even degree part, the invariant in the Jacobi ring $\mathsf{Jac}(\mathbb{A}^5,W)$ of $W$ is spanned by the monomials $1$, $\Pi x_i$, $(\Pi x_i)^2$, and $(\Pi x_i)^3$.
\end{itemize}
Let $h^{q,p}(X;G)$ be the sum of $h^{q,p}(X;g)$ for $g\in G$ and let $\check{h}^{q,p}(X;G)$ be the sum of $\check{h}^{q,p}(X;g)$ for $g\in G$. Let $h^{q,p}([X/G])$ be the sum of $h^{q,p}(X;g)^G$ for $g\in G$ and let $\check{h}^{q,p}([X/G])$ be the sum of $\check{h}^{q,p}(X;g)^G$ for $g\in G$. We put the numbers $h^{q,p}([X/G])$ and $\check{h}^{q,p}([X/G])$ into the form of a diamond. Based on the computations above, the diamonds are
\[
\begin{tikzcd}[row sep=tiny, column sep=tiny]
  &   &     & 1   &     &   &    & &   &   &     &  1  &     &   &   &    \\
  &   & 0   &     & 0   &   &    & &   &   &  0  &     &  0  &   &   &    \\
  & 0 &     & 101 &     & 0 &    & &   & 0 &     &  1  &     & 0 &   &    \\
1 &   & 1   &     & 1   &   & 1, & & 1 &   & 101 &     & 101 &   & 1 &    \\
  & 0 &     & 101 &     & 0 &    & &   & 0 &     &  1  &     & 0 &   &    \\
  &   & 0   &     & 0   &   &    & &   &   &  0  &     &  0  &   &   &    \\
  &   &     & 1   &     &   &    & &   &   &     &  1 , &     &   &   &    
\end{tikzcd}
\]
where $101$ is equal to $1+60\times 1+40\times 1$.

Before we prove Theorem B2, we need the proposition below.

\begin{Proposition}
We put the numbers $\check{h}^{q,p}([X/G])=\dim\HT^{q,p}([X/G])$ into a form of a diamond. Then it is in the form of a Greek cross below
\[
\begin{tikzcd}[row sep=tiny, column sep=tiny]
  &   &      &        & 1     &        &      &    &  \\
  &   &      & 0      &       & 0      &      &    &  \\ 
  &   &\cdots&        & 1     &        &\cdots&    &  \\
  & 0 &      & \cdots &       & \cdots &      & 0  &  \\
* &   & *    &        &\cdots &        & *    &   & * \\
  & 0 &      & \cdots &       & \cdots &      & 0  &   \\
  &   &\cdots&        & 1     &        &\cdots&    &   \\
  &   &      & 0      &       & 0      &      &    &   \\
  &   &      &        & 1.     &        &      &    &
\end{tikzcd}
\]
The vertical line and the horizontal line in the diamond are possibly nonzero and the other part of the diamond is zero. \end{Proposition}
\begin{proof}
We know that the Hodge diamond of a general hypersurface in projective space is in the form of a Greek cross 
%\[
%\begin{tikzcd}[row sep=tiny, column sep=tiny]
%  &   &      &        & 1     &        &      &    &  \\
%  &   &      & 0      &       & 0      &      &    &  \\ 
%  &   &\cdots&        & 1     &        &\cdots&    &  \\
%  & 0 &      & \cdots &       & \cdots &      & 0  &  \\
%* &   & *    &        &\cdots &        & *    &   & * \\
%  & 0 &      & \cdots &       & \cdots &      & 0  &   \\
%  &   &\cdots&        & 1     &        &\cdots&    &   \\
%  &   &      & 0      &       & 0      &      &    &   \\
%  &   &      &        & 1     &        &      &    &
%\end{tikzcd}
%\]
and the numbers (except for the one in the middle degree) in the vertical line are $1$. 

Paragraph~\ref{7.1} shows that the diamond of $\HT^{q,p}([X/G])$ is equal to the Hodge diamond of $X$.
\end{proof}

\paragraph\label{7.5} Similar to the Fermat quintic case, we have an isomorphism of algebras
$$\HH^*(X)\cong(\bigoplus_{g\in\mathbb{Z}/d\mathbb{Z}}\mathsf{Jac}(Y^g,W|_{Y^g})\otimes\omega_g)^{\mathbb{Z}/d\mathbb{Z}}.$$

This isomorphism of algebras is crucial in the rest of this paper because of the reason below. We explained that $\HT^*(X)$ is isomorphic to $\HH^*(X)$ as algebras in the introduction. In the definition of both the simplified and unsimplified products, it is clear that $\HT^*(X)$ is a subalgebra of $\HT^*(X;G)$. Therefore the isomorphism of algebras between the Hochschild cohomology and the Jacobi ring above can help us to study the product on $\HT^*(X;G)$ and on $\HT^*([X/G])$.

Similar to the computation in Paragraph~\ref{diamond}, we conclude that there is an element $\alpha\in \HT^{1,1}([X/G])$ such that $\alpha,\alpha^2,\cdots,\alpha^n$ are all nonzero. This class $\alpha$ corresponds to the monomial $\Pi x_i$ in the Jacobi ring under the isomorphism of algebras above and $\alpha^j$ corresponds to $(\Pi x_i)^j$. Since the numbers in the vertical line of the diamond are $1$ except for the middle of this line, the family $<1,\alpha,\alpha^2,\cdots,\alpha^n>$ determines everything in the vertical line except for the middle of this line.

\paragraph We need to clarify terminology before we prove Theorem B2. When the dimension $n$ of $X$ is an odd number, the horizontal line and the vertical line of the diamond of $\HT^*([X/G])$ do not intersect, and there is nothing special to say. When the dimension $n=2m$ of $X$ is an even number, the vertical line and the horizontal line intersect in the middle $\HT^{m,m}([X/G])$. Let $\alpha$ be the class in $\HT^{1,1}(X)$ in the vertical line. Then the class $\alpha^{m}$ lies in the intersection $\HT^{m,m}([X/G])$ of the vertical line and the horizontal line. The intersection naturally decomposes
$$\displaystyle{\HT^{m,m}([X/G])=\bigoplus_{g\in G} (\HT^{m,m}(X;g)^G).  }$$
Similar to the computation in Paragraph~\ref{diamond}, one concludes that $\HT^{m,m}(X;1)^G$ is one dimensional and it is spanned by $\alpha^m$. Then we have the following
$$\displaystyle{\HT^{m,m}([X/G])=\bigoplus_{1\neq g\in G} (\HT^{m,m}(X;g)^G)\bigoplus<\alpha^m>.  }$$

\begin{Definition}
When the dimension $n$ of $X$ is odd, the horizontal and vertical lines of the diamond do not intersect. We define $VL$ as the vector space spanned by the vertical line and define $HL$ as the vector space spanned by the horizontal line.

When the dimension $n=2m$ of $X$ is even, define $VL$ as the vector space 
\[
\displaystyle{\bigoplus_{i\neq m}\HT^{i,i}([X/G])\bigoplus<\alpha^m>,}
\]
and define $HL$ as the vector space
\[
\displaystyle{\bigoplus_{i\neq0} \HT^{m-i,m+i}([X/G])\bigoplus_{1\neq g\in G} (\HT^{m,m}(X;g)^G).}
\]
\end{Definition}

The intersection $\HT^{m,m}([X/G])$ of the vertical line and the horizontal line is naturally indexed by $g$ \[\HT^{m,m}([X/G])=\bigoplus_{g\in G} (\HT^{m,m}(X;g)^G).
\]
The term $<\alpha>=\HT^{m,m}(X;1)^G$ belongs to $VL$ and all the other terms belong to $HL$.

From the discussion above, we know that the vector space $VL$ is always spanned by $<1,\alpha,\cdots,\alpha^n>$. We also know that $VL$ and $HL$ span the entire diamond of $\HT^{(q,p)}([X/G])$ and that $VL\cap HL=0$.

\begin{Proof}[Proof of Theorem B2]
Let $[X/G]$ be a Calabi-Yau Fermat hypersurface of dimension $n$ and degree $d$. We want to prove that the simplified and unsimplified products on $\HT^*([X/G])$ agree for $[X/G]$.

 Let $*_{s}$ and $*_{u}$ be the simplified product and unsimplified product respectively. Let $\alpha_g$ be a class indexed by $g$ and $\beta_h$ be a class indexed by $h$. To prove Theorem B2, it suffices to show that $\alpha_g*_{s}\beta_h=\alpha_g*_{u}\beta_h$ when $\alpha_g$, $\beta_h$ are both in the vector space $VL$, or both in the vector space $HL$, or one is in $VL$ and the other is in $HL$.
 
When $g$ is the trivial element, the isomorphism of algebras $\HT^*(X)=\HT^*(X,1)\cong\HH^*(X)=\HH^*(X,1)$ is explained in the introduction. One can conclude that $\HT^*(X)^G\cong\HH^*(X)^G$ is a subalgebra of 
$$\displaystyle{\HT^*([X/G])=\bigoplus_{g\in G}\HT^*(X,g)^G }$$
for both the simplified and the unsimplified products using the definition of the products.

It has been shown in Paragraph~\ref{7.5} that the product restricted to 
$$VL=<1,\alpha,\alpha^2,\cdots,\alpha^n>$$
is generated by $\alpha$ as an algebra. Namely, the product restricted to $VL$ is determined by $\HT^*(X)\cong\HH^*(X)$ and it has no contribution from $\HT^*(X;g)$ for nontrivial $g\in G$. Therefore, the simpilifed and unsimplified product agree on $VL$, i.e., $\alpha_g*_{s}\beta_h=\alpha_g*_{u}\beta_h$ when $\alpha_g$, $\beta_h$ are both in $VL$.

When one of $\alpha_g$ and $\beta_h$ is in $VL$ and the other one is in the $HL$, we can assume that $\alpha_g$ is in $VL$, i.e., $g$ is the trivial element $1\in G$ and $\alpha_g=\alpha^i$ for some $i$, without loss of generality. We want to show that the product is either zero or $\alpha_g$ is the unit of this algebra for both the simplified and unsimplified products. The argument for the simplified and the unsimplified products are the same. Let $\cdot$ be the simplified or the unsimplified product in this paragraph. If $\alpha_g$ is not the unit and $\alpha_g\cdot\beta_h$ is not zero, then the degree of $\alpha_g\cdot\beta_h$ is strictly greater than the degree of $\beta_h\in HL$. Therefore $\alpha_g\cdot\beta_h$ is not in $HL$ because all the elements in $HL$ lie in the middle degree. Since the diamond is of the form of a Greek cross, we conclude that $\alpha_g\cdot\beta_h$ must be in $VL$, i.e., it is of the form $\alpha^i$ for some $i$. The class $\alpha_g$ is also of the form $\alpha^i$ for some $i$, so this implies that $\beta_{h}$ is also of the form $\alpha^i$ and it is in $VL$. We obtain a contradiction. This shows that the product is either zero or $\alpha_g$ is the unit of this algebra in this case.
%we can only obtain possibly nonzero product in the following two situations due to degree reason and the fact that the product on $VL$ is completely determined by $\HT^*(X,1)\cong\HH^*(X,1)$. 
%\begin{itemize}
%\item $h$ is also the trivial element.

%\item $h$ is nontrivial element of $G$, and $\alpha_g$ is the unit of the algebra $\HT^*([X/G])$.
    
%\end{itemize}

%In the first case above, $\alpha_g$ and $\beta_h$ lie in %$\HT^*(X,1)\cong\HH^*(X,1)$, so we can conclude that 
%$\alpha_g*_{s}\beta_h=\alpha_g*_{u}\beta_h$. In the second case above, 
Then we look at the case where both $\alpha_g$ and $\beta_h$ are in $HL$. The simplified product must land in $\HH^{n,n}(X)$ and $g$ and $h$ must be inverse to each other because the simplified product preserves the $(q,p)$ bidegree. Recall that there is a class $\gamma_{g,h}$ introduced in the definition of the simplified product. It is the top Chern class of a vector bundle of rank $\iota(g)+\iota(h)-\iota(gh)-\codim(X^{g,h},X^{gh})$.
When $h=g^{-1}$ or $g=1$, the rank above is zero which implies the class $\gamma_{g,h}$ is trivial. As a consequence, we conclude that the class $\gamma_{g,h}$ does not show up in the simplified product on $\HT^*([X/G])$. When the class does not show up in the product, the simplified product is in the form of Proposition~\ref{Prop trivial}. By definition, the unsimplified product must land in
$$\displaystyle{\bigoplus_{i=0}^{c_g}\HT^{n-c_g+i,n+c_g-i}(X),}$$
where $c_g$ is the codimension of $X^g$ in $X$.
However, there is only one nonvanishing term $\HT^{n,n}(X,1)\cong\HH^{n,n}(X,1)$ of degree $2n$ in $\HT^{*}([X/G])$, i.e., the direct summands above is only nonzero when $i=c_g$. Then the unsimplified product must also land in $\HH^{n,n}(X)$. Because of this reason, one concludes that the unsimplified product is also of the form in Paragraph~\ref{Prop trivial}, i.e., we have $\alpha_g*_{s}\beta_h=\alpha_g*_{u}\beta_h$.
\end{Proof}

Before we prove Theorem D, we need the proposition below.

\begin{Proposition}
Let $\alpha_g,\beta_h\in\HT^{*,*}([X/G])$ be a class indexed by $g$ and $h$ respectively. The product $\alpha_g\circ\beta_h$ is determined by the following three cases: both $\alpha_g$ and $\beta_h$ are in $VL$, both of them are in $HL$, and $\alpha_g$ is in $VL$ and $\beta_h$ is in the $HL$. The product restricted to $VL$ is generated by a class $\alpha\in\HT^{1,1}([X/G])$ as an algebra, i.e., it is of the form $<1,\alpha,\alpha^2,\cdots,\alpha^n>$. The product $(\HT^*([X/G]),\circ)$ restricted to $HL$ can be identified with a pairing. This pairing is nondegenerated. When the dimension $n$ is even, this pairing is symmetric, and when the dimension $n$ is odd, this pairing is skew-symmetric. In the last case, the class $\alpha_g$ must be the unit of the algebra or the product vanishes.
\end{Proposition}

\begin{proof}
We consider the first case. Note that the introduced sign $\epsilon(g,h)$ is $1$ in $VL$ because $g$ and $h$ are trivial in this case. Under the identification
$$\HH^*(X)\cong(\bigoplus_{g\in\mathbb{Z}/5\mathbb{Z}}\mathsf{Jac}(Y^g,W|_{Y^g}))^{\mathbb{Z}/5\mathbb{Z}},$$
we know that the space $VL$ is represented by the classes $1$, $\Pi x_i$, $(\Pi x_i)^2$, $\cdots$, and $(\Pi x_i)^n$ in the Jacobi ring, from the explanation above. Moreover, the multiplication structure is preserved~\cite{T}, i.e., the classes are of the form $<1,\alpha,\alpha^2,\cdots,\alpha^n>$, where $\alpha$ is the class represented by $\Pi x_i$.

We consider the second case. The product restricted to $HL$ can be viewed as a pairing. The sign $\epsilon(g,h)$ is introduced to make the product graded commutative. Therefore, it suffices to show that the pairing is nondegenerated. From the previous discussions, we know that the only possibly nonzero product, in this case, is of the form below

$$\HT^{q,p}(X;g)^G\otimes\HT^{q',p'}(X;h)^G\rightarrow H^n(X,\wedge^n T_X)\cong\mathbb{C},$$
where $h=g^{-1}$, $p+p'=q+q'=n$, and $p+q=p'+q'=n$.
%If the product above is always nonzero, then $\HT^{q,p}(X;g)^G$ form a basis for the pairing. Under this basis the pairing is one of the two forms(up to scalars)

%\[\begin{pmatrix}
%    0 & 0 & \cdots & 0 & 1\\
%    0 & 0 & \cdots & 1 & 0\\
%    \vdots &\vdots & 1 &\vdots &\vdots \\
%    0 & 1 & \cdots & 0 & 0\\
%    1 & 0 & \cdots & 0 & 0
%    \end{pmatrix}
%    \]
%    \[
%    \begin{pmatrix}
%    0 & I_m \\
%    -I_m & 0
%    \end{pmatrix}
%    \]
%depending on whether the dimension $n$ is even or odd. Then it is clear that the pairing is nondegenerated. It suffices to prove that the product of the form above is always nonzero.

We expand the term 
$$\HT^{q,p}(X;g)^G=H^{p-\iota(g)}(X^g,\wedge^{q-c_g+\iota(g)}T_{X^g}\otimes\omega_g)^G$$
and similarly the term
$$\HT^{q',p'}(X;h)^G=H^{p'-\iota(h)}(X^h,\wedge^{q'-c_h+\iota(h)}T_{X^h}\otimes\omega_h)^G.$$
The orbifold $[X/G]$ is Calabi-Yau, so we can apply the identification in Theorem C. We have
$$\HT^{q,p}(X;g)^G\cong H^{p-\iota(g)}(X^g,\wedge^{n-q-\iota(g)}\Omega_{X^g})^G,$$

$$\HT^{q',p'}(X;h)^G\cong H^{p'-\iota(h)}(X^h,\wedge^{n-q'-\iota(h)}\Omega_{X^h})^G,$$
and
$$H^n(X,\wedge^n T_X)\cong H^n(X,\wedge^n\Omega_X).$$

Use the definition of the product and the three identifications above to expand the product in detail. The product is the composite map below
$$H^{p-\iota(g)}(X^g,\wedge^{n-q-\iota(g)}\Omega_{X^g})^G\otimes H^{p'-\iota(h)}(X^h,\wedge^{n-q'-\iota(h)}\Omega_{X^h})^G $$
$$\rightarrow H^{p+p'-\iota(g)-\iota(h)}(X^{g,h} \wedge^{n-q-\iota(g)}\Omega_{X^g}\otimes \wedge^{n-q'-\iota(h)}\Omega_{X^h})^G$$
$$\longrightarrow H^n(X,\wedge^n\Omega_X)\cong\mathbb{C}.$$
We are in the case where $h=g^{-1}$. We have $X^{g}=X^{h}=X^{g,h}$ and $\iota(g)+\iota(h)=c_g$. In addition, we have $p+p'=q+q'=p+q=p'+q'=n$. Then the numbers satisfy
$$p+p'-\iota(g)-\iota(h)=n-c_g=d_g,$$
and $$ n+n-q-q'-\iota(g)-\iota(h)=n-c_g=d_g,$$
where $d_g$ is the dimension of $X^g$. The product is greatly simplified thanks to the equations above. The product now can be rewritten as the composite map below
$$H^i(X^g,\wedge^j \Omega_{X^g})^G\otimes H^{i'}(X^g,\wedge^{j'} \Omega_{X^g})^G $$
$$\rightarrow H^{d_g}(X^g,\wedge^j \Omega_{X^g}\otimes\wedge^{j'}\Omega_{X^g})^G\rightarrow H^{d_g}(X^g,\wedge^{d_g} \Omega_{X^g})^G=H^{d_g}(X^g,\wedge^{d_g}\Omega_{X^g})$$
$$\longrightarrow  H^n(X,\wedge^n\Omega_X),$$
where $i+i'=j+j'=d_g$.

The composite map of the first two arrows above is the wedge product. The last arrow is of the form $H^0(X^g,\cO_{X^g})^\vee\rightarrow H^0(X,\cO_{X})^\vee$ by applying Serre duality. Using the definition of the product, one can conclude that the last arrow is exactly due to the natural map 
$$H^0(X,\cO_X)\rightarrow H^0(X^g,\cO_{X^g}).$$
It is clear that the composite map below
$$H^i(X^g,\wedge^j \Omega_{X^g})\otimes H^{i'}(X^g,\wedge^{j'} \Omega_{X^g}) $$
$$\rightarrow H^{d_g}(X^g,\wedge^j \Omega_{X^g}\otimes\wedge^{j'}\Omega_{X^g})\rightarrow H^{d_g}(X^g,\wedge^{d_g} \Omega_{X^g})$$
$$\longrightarrow  H^n(X,\wedge^n\Omega_X)\cong\mathbb{C}$$
is a nondegenerate pairing because it is the standard pairing on de Rham cohomology of each $X^g$. Now we only need to show that it remains nondegenerate after taking $G$-invariants. It follows from the lemma below.

The third case has been proven in the proof of Theorem B2.
\end{proof}

\begin{Lemma}
Let $V$ a vector space and $<,>$ be a nondegenerate pairing on $V$. A finite group $G$ acts on $V$ and preserves the pairing, i.e., $<gv,gw>=<v,w>$ for all $v,w\in W$. Then the induced pairing on the fixed locus $V^G$ remains nondegenerate.
\end{Lemma}
\begin{proof}
Let $w$ be an element in $V^G$. We need to show that if $<w,v>=0$ for all $v\in V^G$, then $w=0$. 

Suppose $w$ is not zero and $<w,v>=0$ for all $v\in V^G$. There exists an element $v_0\in V$ such that $<w,v_0>\neq0$ because the pairing on $V$ is nondegenerated. Then 
$$<w,\frac{1}{|G|}\sum_{g\in G}gv_0>=\frac{1}{|G|}\sum_{g\in G}<w,gv_0>$$
$$=\frac{1}{|G|}\sum_{g\in G}<g^{-1}w,v_0>=\frac{1}{|G|}\sum_{g\in G}<w,v_0>=<w,v_0>\neq0.$$
The first equality in the second row is due to the fact that the $G$ action preserves the pairing. The second equality in the second row is due to the fact that $w$ is invariant under $G$. However, $\frac{1}{|G|}\sum_{g\in G}v_0$ is an element in $V^G$. We get a contradiction.
\end{proof}
Since we have computed the product $(\HT^*([X/G]),\circ)$ explicitly, we are able to prove Theorem D easily as follows.

\begin{proof}[\bf Proof of Theorem D]
As explained at the beginning of this section, there is an isomorphism of algebras

$$\HH^{*}([X/G])\cong \HH^{*}(MF_{gr}(\mathbb{A}^{n+1},G\times \mathbb{Z}/d\mathbb{Z},\Sigma_{i=0}^n \,x_i^d))$$

$$\cong\bigoplus_{g\in G\times \mathbb{Z}/d\mathbb{Z}}(\mathsf{Jac}(Y^g,W|_{Y^g})\otimes\omega_g)^{G\times \mathbb{Z}/d\mathbb{Z}},$$
where $Y$ is $\mathbb{A}^{n+1}$.

To prove Theorem D, it suffices to compare the product $(\HT^*([X/G]),\circ)$ with the product on 
$$\bigoplus_{g\in G\times \mathbb{Z}/d\mathbb{Z}}(\mathsf{Jac}(Y^g,W|_{Y^g})\otimes\omega_g)^{G\times \mathbb{Z}/d\mathbb{Z}}.$$

The product structure on the orbifold matrix factorization category has been defined and studied in~\cite{HLL}.
One can directly check that the products match. In fact, it has been shown in {\em loc. cit.} that the product on the orbifold matrix factorization has a Frobenius algebra structure. Think about the diamond of this algebra which is exactly the same as the diamond of $\HT^{q,p}([X/G])$. Therefore the diamond decomposes into $HL$ and $VL$ as before. The Frobenius algebra structure determines the product structure on $HL$ because it has to be a nondegenerate symmetric or skew-symmetric pairing. There is only one such pairing in a fixed dimension. 

The remaining part of the product can be computed in the same method as the one that has been used in the proof of Theorem B2. The product restricted to $VL$ is again generated by one class $\alpha$ of bidegree $(1,1)$ because the same computation can be done in the Jacobi ring above. Using the same argument in the proof of Theorem B2, one can see that the product of an element $x\in VL$ with an element $y\in HL$ is completely determined due to degree reason. This product $x\cdot y$ has to vanish unless $x$ is the unit of the algebra.
\end{proof}

\paragraph In the rest of this section we explain the possible application of the theorems in this paper to mirror symmetry. We explain that the product on the orbifold polyvector field of $[X/G]$ is expected to match with the product on the state space of the FJRW theory of the Fermat polynomial $W$. 

Standard physical arguments \cite{H03} predict that Fermat hypersurface $X$ is the mirror to $[X/G]$, which in particular implies the topological $A$ model of $X$ is equivalent to the topological $B$ model of $[X/G]$. On the $B$ side we have the complex moduli of $[X/G]$ parameterized by $\psi$ (the defining equation is explicitly given by $\sum x_i^n+\psi\prod x_i=0$), on each point of the moduli there is a 2D TQFT associated with its derived category. By mirror symmetry we would have a corresponding structure on the $A$ side, i.e., we would have a (stringy) Kahler moduli parameterized by $\psi$ and for each $\psi$ a TQFT constructed from symplectic geometry of $X$. This picture is much less understood due to the presence of instanton corrections, we only have precise mathematical definitions in certain limits. The large volume limit of $A$ model corresponds to $\psi\rightarrow\infty$, where the corresponding category is the Fukaya category of $X$ which is closely related to its Gromov-Witten invariants. There is an opposite limit $\psi\rightarrow 0$ which formally means negative infinity Kahler class and the corresponding mathematical objects are the Fukaya-Seidel category and FJRW invariants. In this case of Fermat hypersurface, there is a proof \cite{Ch10} of the expected fact that FJRW theory and Gromov-Witten theory are related by analytic continuation in genus $0$. 

\paragraph In this paper, we studied the $B$ model chiral ring at point $\psi=0$ (note that unlike the large complex structure limit $\psi\rightarrow\infty$, $\psi\rightarrow 0$ converges to a smooth orbifold, categorically this is better behaved) which under mirror symmetry goes to the classical FJRW ring, i.e. we only need the point $0$ instead of its formal neighbourhood. This is purely topological and also expected to be the orbifold Jacobi ring of the mirror Fermat polynomial. Actually, the isomorphism of orbifold Jacobi ring and orbifold Hochschild cohomology has a purely $B$ model explanation without referring to mirror symmetry: for different points in the (stringy) Kahler moduli we get different descriptions of the theory which are canonically isomorphic up to monodromy on the Kahler moduli (this is the mirror to the well-known fact that we need to choose an almost complex structure to define the Fukaya category and different choices are canonically isomorphic). There are two particular points in the Kahler moduli resembling $\psi\rightarrow\infty$ and $\psi\rightarrow 0$ and give two descriptions of the category as (orbifold) derived category and matrix factorization, and their equivalence is the orbifold extension of Orlov's theorem. This equivalence induces an isomorphism of their Hochschild cohomology, which is one of the key ingredients of our proof.

\end{document}